\documentclass[a4paper,11pt]{article}  

\usepackage{graphicx,color,subfigure,multirow, here}
\usepackage{mathtools}
\usepackage[utf8]{inputenc} 
\usepackage[T1]{fontenc}
\usepackage[english]{babel} 
\usepackage{amsmath, amsfonts, amsthm,amssymb,here,dsfont, hyperref, enumitem}
\usepackage{natbib}

\newcommand{\argmin}[1]{\underset{#1}{\operatorname{arg}\!\operatorname{min}}\;}
\newcommand{\one}{\mathds{1}}
 
\newtheorem{theo}{Theorem}[section]

\newtheorem{prop}[theo]{Proposition}

\newtheorem{lemme}[theo]{Lemma}
\newtheorem{rem}[theo]{Remark}
\newtheorem{ass}[theo]{Assumption}

\newcommand{\n}{\noindent }

\newcommand{\E}{\mathbb{E}}
\newcommand{\R}{\mathbb{R}}
\newcommand{\w}{\widehat}

\renewcommand{\P}{\mathbb{P}}

\newcommand{\pen}{\text{pen}}

\makeatletter

\@addtoreset{equation}{section}
\makeatother

\begin{document}

\title{
Nonparametric drift estimation for diffusions\\ with jumps driven by a Hawkes process}
 \maketitle

\begin{center}
Charlotte Dion$^{(1)}$ and Sarah Lemler$^{(2)}$\\
$^{(1)}$ Sorbonne Université, UMR CNRS 8001, LPSM, 75005 Paris, France\\
$^{(2)}$ Laboratoire MICS, \'Ecole CentraleSup\'elec, Université Paris-Saclay 
\end{center}


\begin{abstract}
We consider a $1$-dimensional diffusion process $X$ with jumps. The particularity of this model relies in the jumps which are driven by a multidimensional Hawkes process denoted $N$. 
This article is dedicated to the study of a nonparametric estimator of the drift coefficient of this original process.
We construct estimators based on discrete observations of the process $X$ in a high frequency framework with a large horizon time and on the observations of the process $N$.
The proposed nonparametric estimator is built from a least squares contrast procedure on subspace spanned by trigonometric basis vectors. 
We obtain adaptive results that are comparable with the one obtained in the nonparametric regression context. 
We finally conduct a simulation study in which we first focus on the implementation of the process and then on showing the good behavior of the estimator.
\end{abstract}

\n \textit{Keywords:}
Nonparametric estimator. Model selection. Diffusion. Hawkes process. \\

\n {\bf AMS Classification}: 62G05,  60J60, 60J75 \\

\section{Introduction}\label{sec:intro}

Boosted by applications in neuroscience and finance this paper is dedicated to a jump-diffusion model.
We focus on a new kind of diffusion with jumps driven by a multidimensional Hawkes process. 
The process $X$ of interest is solution of the following equation
\begin{eqnarray}\label{eq:model}
 dX_t =b(X_t)  dt + \sigma(X_t)  dW_{t} + a(X_{t^-}) \sum_{j =1}^M dN^{(j)} _t,
\end{eqnarray}
where $X_{t^-}$ denotes the process of left limits, $N=(N^{(1)}, \ldots, N^{(M)})$ is a $M$-dimensional Hawkes process with (conditional) intensity function $\lambda$ (that will be detailed in the next section) and $W$ is the standard Brownian motion independent of $N$. This process is a diffusion process between the jumps of $N$. 

The first probabilistic properties of this process have been established in \cite{DLL}. The present paper is part of the continuation of this first work and focuses on the estimation part. Since this model is new, no statistical inference results have yet been carried out on it. In this paper, we focus on the nonparametric estimation of the drift function which leads the dynamic of the trajectories and which is already a real challenge. An estimator obtained from a minimization of a modified least squares contrast is proposed. The difficulty is to control the variance of the estimators taking into account the jumps in order to establish an oracle-type inequality for the empirical risk. 

\subsection{Motivation}

A main motivation is to find a good model to describe neurons membrane potentials between spikes. A spike is defined as an electrical signal that exceeds a certain threshold. Considering one neuron, jump diffusion models are often used to fit the dynamic of the membrane potential of a neuron \citep[see \emph{e.g.}][]{JBHD}. Indeed, according to the literature \citep[see \emph{e.g.}][]{HOP} diffusion processes seem suitable to model the membrane potential of a neuron between successive spikes. Then, Poissonian firing can be added (with a state dependent firing rate) to the diffusion process in order to describe the whole dynamic of the neuron including its own spikes. In the present work, we keep the diffusion model to describe the segments of the membrane potential between spikes of a fixed neuron, but, enriched by a Hawkes jump term, representing the influence of a network of neurons surrounding the one neuron of interest in the study. In fact, we take advantage of both kind of data available due to the biologists. On one hand, we can recorded the membrane potential of one neuron at high frequency over a long time interval (intracellular recording), this signal is continuous. On the other hand, we are also able to record the spike trains occurence of $M$ neurons (extracellular recording), which are the sequence of times when the $M$ neurons spike; we obtain a discrete signal.  The network activity is captured in our model \eqref{eq:model} by the multidimensional Hawkes process. This hybrid model, in which the dynamic of the membrane potential jumps when one of the $M$ neurons around spikes,
allows to combine the wealth of both kind of data by using them all together in a single model.

But we also believe in the usefulness of this new model to describe many other phenomena. In financial mathematics for example,
this model could enrich the collection of models for the stochastic volatility of an asset price. The asset price along time $S_t$ can be modeled by a stochastic differential equation and the volatility coefficient $\sigma_t$ \citep[see \emph{e.g.}][]{TOUZI} can be written as $f(X_t)$ where $X_t$ satisfies Equation \eqref{eq:model}. In this case, the arrival of economic news are well described by the Hawkes process and affects the dynamic of asset prices  \citep[see \emph{e.g.}][]{RAMBALDI}.

The idea here is different from interaction diffusion models \citep[see][]{CF13,GOBET} because it takes into account a continuous signal and a discrete signal all together. For what concerns for example the neuronal application, the data are the jump processes observed from several neurons and the membrane potential of one neuron.

The main contribution of this model is that the Hawkes process has a special time dependence structure:
the conditional intensity $t \mapsto \lambda_t$ is a random predictable function that depends on the past before time $t$.
Nevertheless, we do not assume that the jump intensity of $N$ at time $t$ depends on the state of the process $X_t$.

\subsection{State of the art}

Diffusions find wide use in the modeling of dynamic phenomenons in continuous time.
Diffusion processes (without jumps) are Markovian processes, used in applied problems as medical sciences \citep[see \emph{e.g.}][]{DS2013}, physics  \citep[see \emph{e.g.}][]{PHY} and financial mathematics  \citep[see \emph{e.g.}][]{vasicek, ELK}. Then, jump-diffusion processes are born for example in the risk management context \citep[see][]{TANKOV} driven by a Poisson process or a L\'evy process \citep[][]{ MASUDA, mancini, SCHMISSER2}. 
In neuroscience the trajectories of the membrane potential of one fixed neuron have been modeled through different diffusion processes \citep[][]{HOP,JBHD, DION2014}. Many other models exist, taking into account other quantities additionally to the membrane potential \citep[for example the Hodgkin-Huxley model in the recent work][]{HLT}. 

The Hawkes process \citep[][]{HAWKES71} comes later analogically to auto-regressive model. It has first been introduced to model earthquakes and their aftershocks in seismology \citep{Vere82}. This  point process generalizes the Poisson process. It has been used in genomics \citep{bonnet} and  more recently it has been considered to model the intensity of spike trains of neuronal networks \citep{RBRTM, DL2016}. In this last context the interaction function are often piecewise constant, see e.g. \citep{BM1996}.

In this last case, it is a  multidimensional Hawkes also called mutually exciting point processes that is considered to model the interactions between several neurons \citep[see][for theoretical studies]{BM1996, DFH}.
This model can deal with observations which represent events associated to agents or nodes on a given network, and that arrives randomly through time but that are not stochastically independent. 
There is a large literature in the financial side on Hawkes processes, among others recently \cite{RJ2015, BACRYFINANCE}. 
Moreover, social networks are nowadays often modeled using a Hawkes process \citep[see \textit{e.g.}][]{twit}.
The strength of the multidimensional Hawkes process is mainly the time dependency of each process, which is called sometimes "self-excitation" and modeled by the kernel function.
The inference of this process has been investigated parametrically for example in \cite{BACRYFINANCE}, nonparametrically in \cite{BACRYkernel,LV,kirchner2017estimation} and  through neighbor method in \cite{Hansen15}  and \cite{DGL2016}.

In this paper we focus on the case of linear exponential Hawkes process. It has several advantageous properties. First theoretically this choice implies the Markov property for the intensity process of $N$. Secondly,
 the expected value of arbitrary functions of the point process is explicit
\citep[see][]{BACRYFINANCE} and it can be exactly simulated \citep[][]{DASSIOS}. 
The allows a clustering representation of the Hawkes process \citep[see][]{HO1974,RRB} and non-negativity of the kernel functions eliminates the inhibition behavior. 
Besides, it is in this context that ergodicity and $\beta$-mixing properties have been proved for $X$ in \cite{DLL}. The following results are based on these characteristics.

\subsection{Main contribution}

Statistical inference for diffusions with jumps contains many challenges. Model \eqref{eq:model} is presented in detail the next section. It is decomposed in three terms: a drift term, a diffusion term and a jump term. The jump term is contained in the multidimensional Hawkes process. The model have continuous-time dynamic but it has to be inferred from discrete-time observations (whether one or more jumps are likely to have occurred between any two consecutive observation times). We assume to observe the process $X$ at equidistant times with time step $\Delta$ which goes to zero when the number of observations $X_{k\Delta}$ named $n$ goes to infinity. Besides, the horizon time $T=n \Delta \rightarrow \infty$ with $n$. We also assume to observe the jumps times of $N$ on $[0,T]$.

The main purpose of this paper is to propose a nonparametric estimator for the drift function $b$ in this new framework.
There is a few frequentist nonparametric work where the drift function of jump-diffusion processes is estimated \citep[][]{SY,H99, mancini, SCHMISSER3, amorino2018contrast}, but to the best of our knowledge it is the first time that jumps are driven by a Hawkes process. Nevertheless, it is a well known problem for diffusion processes: one can cite for example \cite{bibby} and \cite{kesslersorensen} for martingale estimation functions, \cite{gobetetal} in the low frequency context,  \cite{H99,CGCR} for least squares contrast estimator.

To make our input in this statistical field, we focus on the nonparametric estimation of the drift function named $b$,
when the coupled process $(X, \lambda)$ is ergodic, stationary and exponentially $\beta$-mixing. The proposed method is based on a simple penalized mean square approach. The particularity relies in the contrast function which takes into account the difficulty of the additional Hawkes process. The estimators are chosen to belong to finite-dimensional spaces of variable dimension.
An upper bound for the empirical risk of the collection of estimator is obtained. Then, we deal with an adapted selection strategy for the dimension. The final estimator achieves for example the classical speeds of convergence for Besov regularity (without logarithmic loss) under some assumptions.

\subsection{Plan of the paper}

First, the model is introduced in Section \ref{sec:notations}. Then, the estimator is constructed and studied along Section \ref{sec:estimation}. 
A numerical study is  presented in Section \ref{sec:simu}. We lead a discussion in Section \ref{sec:discussion}.
Lastly, the proofs are detailed in Section \ref{sec:proof}. In this last section some more technical lemmas are exposed.

\section{Framework and assumptions}\label{sec:notations}

\subsection{The Hawkes process}

Let $(\Omega,\mathcal{F},\P)$ be a probability space. 
Let us define the 
Hawkes process on $\R^{+}$-time via conditional stochastic intensity representation.
We denote 
the $M$-dimensional point process $N_t=(N_t^{(1)},\dots,N_t^{(M)})$. Its intensity is a vector of non-negative stochastic intensity functions given by
a collection of baseline intensities, which are positive constants  $ \xi_j , 1 \leq j \leq M ,$ and $M\times M$ interaction functions $ h_{i, j } : \R ^+ \to \R^+ , 1 \le i, j \le M ,$ which are measurable functions. 
Let moreover $n^{(i)} , 1 \le i \le M,$  be discrete point measures on $\R^{-} $ satisfying that 
$$ \int_{\R^{-} } h_{i,j}  ( t- s) n^{(i)} (ds ) < \infty \mbox{ for all } t \geq 0 .$$
We interpret them as {\it initial condition} of our process. 
The linear Hawkes process with parameters $ ( \xi_i , h_{i j } )_{ 1 \le i , j \le M} $ and with initial condition $ n^{(i)}, 1 \le i \le M,$ is a multivariate counting process $N_t ,~ t \geq 0,$ such that 
$\P-$almost surely, for all $ i \neq j  ,  N^{(i)} $ and $ N^{(j) }  $ never jump simultaneously,
and 
for all $1 \le i \le M, $ the compensator of $N^{(i) }_t $ is given by $\Lambda^{(j)}_t:=\int_0^t \lambda^{( i) }_s ds  $ where 
\begin{equation}\label{eq:intensity}
\lambda^{(j)}_t= \xi_j + \sum_{i=1}^M \int_{0}^{t-} h_{i,j}(t-u)dN^{(i) }_u  + \sum_{i=1}^M \int_{-\infty}^{0} h_{ i,j }(t-u)dn^{(i) }_u.
\end{equation}
Here $N_t^{(j)}$ is the cumulative number of events in the $j-$th component at time $t$ and $dN^{(j)}_t$ represents the number of points in $[t, t+dt]$.
The associated measure can also be written
$N^{(j)}(dt)=\sum_{k \in \mathbb{N}} \delta_{T^{(j)}_k}(dt)$.
This intensity process $\lambda$ of the counting process $N$ is the $\bar{\mathcal{F}}_t-$ predictable process such that $(\widetilde{N}_t=N_t-  \Lambda_t)_t$ is a $\bar{\mathcal{F}}_t-$local martingale with $\bar{\mathcal{F}}_t=\sigma\{N_s, 0\leq s\leq t\} =\bigvee_{j=1}^M {\mathcal{F}}_t^{(j)}$ the history of the process  $N$ \citep[see][]{DVJ}.

If the functions $ h_{i,j } $ are locally integrable, the existence of a process $(N^{(j)}_t)_{t \geq 0}$ with prescribed intensity  on finite time intervals follows from standard arguments \citep[see \emph{e.g.}][]{DFH}.
Then, $\xi_j$ is the exogenous intensity and we denote $(T_k^{(j)})_{k\geq 1}$ the non-decreasing jump times of the process 
$N^{(j)}$.

The interaction functions $h_{i, j}$ represent the influence of the past activity of subject $i$ on the subject $j$, the parameter $\xi_j$ is the spontaneous rate and is used to take into account all the unobserved signals. 
In the following we focus on the exponential kernel functions defined by 
\begin{equation}\label{eq:defh}
h_{i, j}: \R^+ \rightarrow \R^+, ~h_{i,j}(t)= c_{i,j}e^{-\alpha t}, ~\alpha >0, ~c_{i,j}>0, ~\quad 1 \leq i, j \leq  M.
\end{equation}
%
The conditionnal intensity process $(\lambda_t)$ is then Markovian. Before the first occurrence and when $n^{(i)}\equiv 0$ in \eqref{eq:intensity}, which means that the Hawkes process has empty history,, all the
point processes $N^{(j)}$ behave like homogeneous Poisson processes with constant intensity $\xi_j$. But as soon as an occurrence
appears for a process $N^{(i)}$, then this affects all the process by increasing  the conditional
intensity via the interaction functions.

\subsection{General model assumptions}

We are now able to write the process as $M+1$ stochastic equations
$$
\begin{cases}
d\lambda^{(j)}_t &= -\alpha ( \lambda^{(j)}_t- \xi_j) dt + \sum_{i =1}^M c_{i, j}dN^{(i)}_t,   \quad j=1, \dots, M \\
dX_t&= b(X_t) dt+ \sigma(X_t)dW_t +  a(X_{t^-})\sum_{j=1}^Md{N}^{(j)}_t
\end{cases}
$$
with $\lambda_0^{(j)}, ~ X_{0}$ random variables independent from the other variables, then we see that 
$(\lambda^{(1)}_t,\dots, \lambda_t^{(M)}, X_t)$ is a Markovian process for the general filtration
$$ \mathcal{F}_t= \sigma( \xi_j, W_s,{N}_s^{(j)},j=1, \dots M, 0\leq s \leq t).
$$

\n The  process is observed at high frequency on the time interval $[0,T]$, the observations are denoted $X_0, X_{\Delta}, \dots X_{n\Delta}$, 
with $\Delta \rightarrow0, n \rightarrow \infty \text{ and } T:= n \Delta \rightarrow \infty$, and the jumps times $T_k^{(j)}, j=1, \dots, M$ on $[0,T]$.

The size parameter $M$ is fixed and finite all along  and asymptotic properties are
obtained when $T \rightarrow \infty$.

\begin{ass}\label{ass:coeff}
Assumptions on the coefficients of $X$.
\begin{enumerate}
\item $a,b, \sigma$ are globally Lipschitz, and $ b $ and $ \sigma $ are of class $\mathcal{C}^2. $
\item There exist positive constants $c, q $ such that for all $x \in \R, $ $  | b' ( x) | +| \sigma' (x) | \le c  $ and  $ | b'' ( x) | + | \sigma '' (x)| \le c ( 1 + |x|^q )  .$    
\item There exist positive constants $a_1 $ and $\sigma_0, ~\sigma_1, $ such that $a(x) < a_1$ and $0< \sigma_0< \sigma(x)<\sigma_1$ for all $ x \in \R.$ 
\item There exist $d \geq 0$, $ r > 0 $ such that for all  $x$ satisfying that $ |x| >r$, we have $xb(x) \leq -d x^2.$
\end{enumerate}
\end{ass}
\n Under the three first assumptions, Equation \eqref{eq:model} admits a unique strong solution
(the proof can be adapted from \citep[][]{LEGALL} under the Lipschitz assumption on function $a$). 
The fourth additional assumption
is classical in the study on the longtime behavior of $X$ \citep{Has, veret} and ensure its ergodicity \citep[see][]{DLL}.

\begin{ass}\label{ass:h}
Assumptions on the kernels.
\begin{enumerate}
\item 
Let $H$ be the matrix with entries $H_{i,j}= \int_0^\infty h_{i,j}(t)dt=c_{i,j}/\alpha$, $1 \leq i,j, \leq M$. The matrix $H$ has a spectral radius smaller than $1$. 
\item The offspring matrix $H$ is invertible. 
\item For all $i,j\in\{1,\dots,M\}$, $\alpha\geq c_{i,j}$. 
\end{enumerate}
\end{ass}
Assumption \ref{ass:h}.1. implies that 
$(N_t)$ admits a
version with stationary increments \citep[see][for e.g.]{BM1996}. In the following we always consider that this assumption holds.
Also we consider that $(N_t)$ correspond to the asymptotic limit and $(\lambda_t)$ is a stationary process.
The second Assumption \ref{ass:h}.2. is required to ensure the positive Harris recurrence of the coupled process $(X_t, \lambda_t)$ discussed in the next section.
The last Assumption is needed to prove Proposition \ref{resultatcle}.
For example, under the condition $( \sum_{i,j} c_{i,j}^2)^{1/2} < \alpha$ Assumptions \ref{ass:h}.1. and 3. hold.

Besides, as shown in \cite{BDH}, under Assumption \ref{ass:h} 
$$\E[\lambda_t]= (I_M- H)^{-1} \xi, 
$$
where $I_M$ the identity matrix of size $M$, $\xi= ^t(\xi_1, \dots, \xi_M)$, 
and for $\delta >0$, 
\begin{equation}\label{eq:bacry}
\E[N_{t+\delta}-N_t]= \delta(I_M -H)^{-1} \xi.
\end{equation}

\subsection{First results on the model}

We refer to the paper \cite{DLL} for the proofs of the following probabilistic results.

\begin{theo}[Dion, Lemler, L\"ocherbach (2019)]\label{thm:ergo}
Under Assumptions \ref{ass:coeff}, \ref{ass:h} the process $ (X_t, \lambda_t)_{t \geq 0} $ is positive Harris recurrent with unique invariant measure $ \pi$. In particular, for any starting point $(x,y) $ and any positive measurable function $g : \R \times \R^M \to \R_+, $  as $T \rightarrow \infty, $ $\P_{(x,y)}-$almost surely,
$$ \frac{1}{T} \int_0^T g(X_t, \lambda_t)dt\to  \pi(g),$$
where 
$ \displaystyle
\pi(g)=\int g(u,v)\pi(du,dv). 
$
\end{theo}
Moreover, in \cite{DLL}, the proven 
 Foster-Lyapunov type condition implies that for all $t \geq 0$, $\E[X_t^4]< \infty$.

\begin{ass}\label{ass:statio}
$(X_0,\lambda_0)$ has distribution $\pi$.
\end{ass}
\n The process $(X_t, \lambda_t)$ is then in its stationary regime in the following. Note that $(X_t)$ is non Markovian but invariant anyway.\\

Recall that the $\beta-$mixing coefficient of the stationary Markovian process $(Z_t)=(X_t, \lambda_t)$ is given by 
$$  \beta_{Z}(t)=  \int \| P_t ( z, \cdot ) - \pi  \|_{TV} \pi ( dz )$$
where 
$\displaystyle \| \mu\|_{TV} := \sup_{ g : |g | \le 1 } \mu ( g ) $.
Moreover,
$ \beta_X (t) = \int_{\R \times \R^M}  \| P^1_t ( z, \cdot ) - \pi^X  \|_{TV} \pi ( dz )$ where
$P_t^1 (z, \cdot ) $ is the projection on $X$ of $P_t (z, \cdot )$, meaning $P_t^1 (z, dx ) = P_t (z, dx \times \R^M)$, similarly  $\pi^X (dx) = \pi ( dx \times \R^M)$ is the projection on the coordinate $X$. 
Then, 
under the previous Assumptions \ref{ass:coeff},  \ref{ass:h} and \ref{ass:statio}, according to \cite{DLL} Theorem 4.9 the process $(Z_t)=(X_t,\lambda_t)$ is exponentially $\beta-$mixing and there exist some constants $K, \theta > 0 $ such that 
\begin{equation}\label{eq:betamix}
\beta_X (t) \leq  \beta_{Z}(t)\leq Ke^{- \theta t } .
\end{equation}
%

The following result is analogous to Proposition A obtained in \cite{GLOTER}. It is very useful for the following estimation part.
\begin{prop}\label{resultatcle}
Grant Assumptions \ref{ass:coeff}, \ref{ass:h} and \ref{ass:statio}. For all $t \geq 0$, $0<\delta<1$, there exists a constant $C(M, a_1)>0$ such that for $p=2$ or $p=4$ we have
$$
\E\left[\underset{s \in [t, t+\delta]}{\sup}\left| X_s-X_t\right|^{p} \right] \leq C(M, a_1)\delta.
$$
\end{prop}
\n The proofs are relegated in Section \ref{sec:proof}. \\

The following section is dedicated to the estimation of the function $b$.
We assume that functions $\sigma, a$ are known together with the constants 
$\xi_j,~ \alpha, ~c_{i,j}, ~  (i,j) \in \{1,\dots, M\}$, and the function $b$ is unknown.

 \section{Estimation procedure of the drift function}\label{sec:estimation}

The functions $b$ is estimated only on a compact set $A$ of $\R$. 
It follows from Proposition 4.7 of \cite{DLL} that $\pi^X$ the projection of the invariant measure $\pi$ onto the $X-$ coordinate is bounded from below on any compact subset of $\R$, and we denote: 
\begin{equation}\label{eq:bornepi}
0<\pi_0 \leq \pi^X(x).
\end{equation}
Note that $\pi^X$ is also bounded from above if $a \equiv 0$ or if $h_{i,j} \equiv 0$ for all $i,j$.
Here if $b$ is bounded then the result also holds (see details in \cite{DLL}, it comes from the bounds of the transition densities given in \cite{GOBETIHP}).
Moreover, we denote: 
$ \|\cdot\|_{\pi^X} $ the $\mathbb{L}^2$-norm of $\mathbb{L}^2(A,\pi^X(x)dx)$.

\subsection{Non-adaptive estimator and space of approximation}

From the observations, the following increments are available:
\begin{equation}\label{eq:Y}
Y_{k\Delta}:= \frac{X_{(k+1)\Delta}-X_{k\Delta}}{\Delta}=b(X_{k\Delta})+ I_{k\Delta}+ {Z_{k\Delta} }+ T_{M,k\Delta}
\end{equation}
with 
\begin{eqnarray}\label{IZT}
I_{k\Delta}&:=&  \frac{1}{\Delta} \int_{k\Delta}^{(k+1)\Delta} (b(X_s)-b(X_{k\Delta})) ds, \quad Z_{k\Delta}:=\dfrac{1}{\Delta}\int_{k\Delta}^{(k+1)\Delta} \sigma(X_s)dW_s,\\
 T_{M, k\Delta}&:=& \frac{1}{\Delta} \int_{k\Delta}^{(k+1)\Delta}a(X_{s-})  \sum_{j=1}^M dN_s^{(j)}.
 \end{eqnarray}
 Considering the Doob-Meyer decomposition, the $T_{M,k\Delta}$ term can be decomposed in 
$$\displaystyle 
T_{M, k\Delta}=  \frac{1}{\Delta} \int_{k\Delta}^{(k+1)\Delta}a(X_{s-}) \sum_{j=1}^M {\lambda}_s^{(j)}ds+ \frac{1}{\Delta} \int_{k\Delta}^{(k+1)\Delta} a(X_{s-})\sum_{j=1}^M d\widetilde{N}^{(j)}_s.$$
Following \cite{Hansen15} the second term of this decomposition can be seen as a "noise" term.
In decomposition \eqref{eq:Y} the first term is the term of interest. The second term is negligible, the third term is a centered noise term. But the fourth term $T_{M, k\Delta}$ is then not centered and not negligible. The real regression-type equation is actually:
$$Y_{k\Delta}-T_{M, k\Delta}=b(X_{k\Delta})+I_{k\Delta}+ Z_{k\Delta}.$$

Based on these variables, we propose a nonparametric estimation procedure for the drift function $b$ on a closed interval $A$. 
We consider $\mathcal{S}_m$ a linear subspace of $\mathbb{L}^2(A)$ such that $\mathcal{S}_m= \text{span}(\varphi_1, \dots, \varphi_{D_m})$ of dimension $D_m$ 
where $(\varphi_\ell)_\ell$ is an orthonormal basis of $\mathbb{L}^2(A)$.
We denote $\mathcal{S}_n=\underset{m\in\mathcal{M}_n}\bigcup \mathcal{S}_m$ where $\mathcal{M}_n$ is a set of indexes for the model collection. The contrast function is defined for $t \in \mathcal{S}_n$, by
\begin{eqnarray}\label{eq:contrast2}
\gamma_{n,M}(t)&:=&
\frac{1}{n} \sum_{k=1}^{n} (U_{k\Delta}-t(X_{k\Delta}))^2 , \quad U_{k\Delta}:=Y_{k\Delta}-T_{M, k\Delta}.
\end{eqnarray}
The associated mean squares contrast estimator is 
\begin{eqnarray}\label{estimb}
 \w{b}_m &:=& \argmin{t\in \mathcal{S}_m}{\gamma_{n,M}(t)}.
\end{eqnarray}
Let us comment this nonparametric mean squares contrast $\gamma_{n,M}$. In the literature of nonparametric frequentist drift estimation from discrete data \citep{H99, CGCR} the regression is done on the random variables $Y_{k\Delta}$'s. We consider rather a regression on $U_{k\Delta}$ which depends on $T_{M,k\Delta}$ to take into account the jumps of the process. This new contrast leads to a better approximation of the coefficient $b$ according to formula \eqref{eq:Y} when the coefficient $a$ is known. 
For example in the neuronal context presented in Introduction, $a$ can be assumed to be linear at first to represent that the impact of the neuronal network surrounding the neuron of interest is linear in the position of it. Then, as the jump process is observed and  also $X$ at discrete times, the coefficient $T_{M, k\Delta}$ can be approximated.

\begin{rem}
The estimator $\w{b}_m$ is defined by $\displaystyle \sum_{\ell=1}^{D_m} \w{\alpha}_\ell \varphi_\ell$ and $\alpha$ is the solution of the equation
 $^t \Phi_m U= (^t \Phi_m \Phi_m)\alpha$ where $\Phi_m$ is the Gram matrix $\Phi_m=(\varphi_j(X_i))_{i,j})_{i=1, \ldots, n; ~j=1, \ldots, D_m}$ of size $n \times D_m$ and $U=(U_\Delta,U_{2\Delta}, \dots,U_{n\Delta})$.
Thus, the solution may be not unique if the kernel of $\Phi_m$ is not reduced to $\{0\}$. Nevertheless, the sequence $(\w{b}_m(X_{\Delta}), \dots \w{b}_m(X_{n\Delta}))$ is unique and it is this sequence that we take into account in the following empirical risk .
\end{rem}

%
The spaces of approximation $\mathcal{S}_m$ must satisfy the following  key properties.
\begin{ass}\label{ass:basis}
Assumptions on the subspaces.
\begin{enumerate}
\item $\exists \phi_1>0, ~ \forall t \in \mathcal{S}_m, ~ \|t\|_\infty ^2\leq \phi_1 D_m \|t\|^2$
\item The spaces $\mathcal{S}_m$ are nested spaces of  dimension $D_m \leq D_n$ such that there exists $\mathcal{S}_n$ of dimension $D_n$ such that $\mathcal{S}_m \subset\mathcal{S}_{m+1} \subset \mathcal{S}_n$ for all $m \in \mathcal{M}_n$.
\item There exists a finite positive constant $\phi_{\pi^X}$ such that $$\phi_{\pi^X}:= \sup_{j=1, \dots, D_n} \int_A\varphi_{j}^2(x) \pi^X(x) dx = \sup_{j=1, \ldots, D_n} \| \varphi_j\|^2_{\pi^X} .$$
\end{enumerate}
\end{ass}

Note that Assumption \ref{ass:basis}.3 is automatically true for example for the trigonometric basis (also called Fourier basis) as we have $\E[X_t^2] <\infty$.

\subsection{Risk bound}

The next paragraph is dedicated to the main result obtained on the non-adaptive estimator of the drift function on $A$: $\w{b}_m$ for all model $m$. 
We consider the empirical risk $\E[\|t-b\|_n ^2]$ define through the following empirical squared norm:
$$\|t\|^2_n: = \frac{1}{n} \sum_{k=1}^{n} t^2(X_{k\Delta}).$$
\begin{prop}\label{prop:MISE}
Under Assumptions \ref{ass:coeff}, \ref{ass:h}, \ref{ass:statio} and \ref{ass:basis}, if $\Delta \rightarrow 0$ and $(n\Delta)/(\ln^2(n)) \rightarrow \infty$ when $n \rightarrow \infty $, if $D_n \leq O(\sqrt{n\Delta}/\ln(n))$, then,
 the estimator $\w{b}_m$ of $b$ on $A$ given by Equation \eqref{estimb} satisfies
$$\E[\| \w{b}_m-b\|_n ^2] \leq 13 \inf_{t \in \mathcal{S}_m}\|b-t \|_{\pi^X}^2+C_1 \frac{D_m}{n\Delta}+C_2\Delta + \frac{C_3}{n\Delta}$$
with $C_1$ depending on $\phi_1, \sigma_1$, and $C_2$ depending on $\sigma_1$, $C_3$ depending on $a_1, \sigma_1, M$.
\end{prop}
Let us precise here that the notation $b$ is for $b\one_A$ as $b_m$ and $\w{b}_m$ are $A$-supported.
The upper bound of Proposition \ref{prop:MISE} is decomposed into different types of error. The first term is the bias term which naturally decreases with the dimension of the space of approximation $D_m$. The second term is the variance term or the estimation error which increases with $D_m$. The third and fourth terms come from the  discretization error.

Moreover, if $\pi^X$, is bounded from above by $\pi_1$ on $A$, then we can provide  the classical nonparametric rate of convergence depending only on $n,\Delta$. Indeed,
let us assume that $b$ belongs to a Besov ball denote $\mathcal{B}_{\alpha, 2, \infty}(A)$ \citep[see the reference book][]{DL93}, where $\alpha$ measures the regularity of the function. 
Then, taking $t= b_m$ the projection of $b$ on $\mathcal{S}_m$, $\|b-b_m\|^2 \leq C(\alpha, L)D_m^{-2\alpha}$, and choosing $D_m=(n\Delta)^{1/(2\alpha+1)}$  with  $\Delta=o(1/(n\Delta))$ leads to 
$$\displaystyle \E[\| \w{b}_m-b\|_n ^2]  = O ((n\Delta)^{-2\alpha/(2\alpha+1)}).$$
This is the optimal nonparametric rate of convergence in the case of estimation of the drift without the additional jumps coming from the multivariate Hawkes process \citep[see][]{H99}.
\\

\begin{rem}
Besides, let us emphasized that the result is expressed for the expectation of the empirical norm because of the difficulty of the problem. Indeed, as only one observation of the process is available on $[0,T]$, it is tricky to obtain theoretical guaranties for the expectation of the integrated norm. Indeed, we are not exactly in the regression context because the observations $(X_{k\Delta},U_{k\Delta})_k$ are dependent. Thus, the consistency results are obtained for the empirical norm as in \cite{BCV, BCV2, CGCR}. 
 While, in the (independent) regression context, a truncated version of the estimator is convergent of the integrated norm (see \cite{baraud} or \cite{GYORFI} Chapter 12). 
 
Nevertheless, if we allow ourselves to make repeated observations of the process( which is beyond the scope of the paper), then, shorting the collection of models, might provide the consistency of the estimator for an integrated norm \citep[see the recent work][]{CGC2019}.
\end{rem}

The objective of the following section is to obtain a oracle-type inequality for the empirical risk for a final estimator where the dimension is chosen through a data-driven procedure.

\subsection{Adaptation procedure}

We must define a criterion in order to select automatically the best dimension $D_m$ (the best model) in the sense of the empirical risk. This procedure should be adaptive, meaning independent of $b$ and dependent only on the observations. 
The final chosen model minimizes the following criterion
\begin{equation}\label{eq:estimadapt}
\w{m}:= \argmin{m \in \mathcal{M}_n}{\left\{ \gamma_{n,M}(\w{b}_m)+ \text{pen}(m) \right\}}
\end{equation}
with $\pen(\cdot)$ the increasing function on $D_m$ given by 
\begin{equation}\label{eq:pen}
\pen(m): = \rho \sigma_1^2\frac{ D_m}{n\Delta},
\end{equation}
where $\rho$ is a universal constant that has to be calibrated for the problem (see Section \ref{sec:simu}).

\begin{theo}\label{theo:adapt}
Under Assumptions \ref{ass:coeff}, \ref{ass:h}, \ref{ass:statio} and \ref{ass:basis}, if $\Delta \rightarrow 0$ and $(n\Delta)/(\ln^2(n)) \rightarrow \infty$ when $n \rightarrow \infty $, if $D_n \leq O(\sqrt{n\Delta}/\ln(n))$, then, 
for the final estimator $\w{b}_{\w{m}}$ of $b$ on $A$ given by \eqref{estimb} and \eqref{eq:estimadapt},
there exists $\rho_0$ such that for $\rho\geq 7\rho_0$, 
  $$\E[\|\w{b}_{\w{m}}-b\|_{n}^2] \leq  
C\underset{m\in \mathcal{M}_n}{\inf} \left\{ \inf_{t \in \mathcal{S}_m}\|b-t\|_{\pi^X}^2+\pen(m)\right\}+ C'\left(\frac{1}{\Delta n}  +\Delta\right).
  $$
where ${\rm pen(.)}$ is given by (\ref{eq:pen}), and $C$ is a numerical constant independent from $n$ and $D_m$ and $C'$ depends on the parameters of the model.
\end{theo}

The result teaches us that the estimator $\w{b}_{\w{m}}$ realizes automatically the best compromise between the bias term and the penalty term which has the same order than the variance term.
Again a similar result is obtained for diffusion processes in \cite{CGCR} which confirms that when the jump process is observed the estimation of the drift function is not modified.

\section{Numerical study}\label{sec:simu}
%
%
%

In this section we describe the protocol of simulations which has been used to evaluate the presented methodology.

\subsection{Implementation of the process and chosen examples}

\begin{figure}
\centering
\includegraphics[height=7cm, width=14cm]{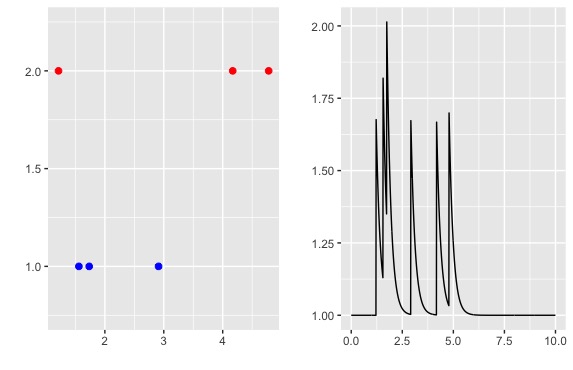}
\caption{Example of simulation of the 2-dimensional Hawkes process: left the two trains of jump times, right the sum of the intensities $(\lambda^{(j)}_t)_t$ on $[0;10]$.}
\label{fig:hawkes}
\end{figure}

We simulate the Hawkes process $N$ with $M=2$ and here we denote $(T_k)_k$ the sequence of jump times of the aggregate process (all the jump times of $N$ sorted and gathered). 
For this simulation we use a Thinning method.
One could also use the \textbf{R}-package 
 \texttt{hawkes} for classical kernels or the package \texttt{tick} for Python language. Also, an exact simulation procedure is possible for the exponential Hawkes \citep[see \emph{e.g.}][]{DASSIOS}.
The intensity process is written as
$$\displaystyle \lambda_t^{(j)}= \xi_j 
 +(\lambda_0^{(j)} -\xi_j) e^{-\alpha t}
+ \sum_{T_k^{(j)} < t} c_{j,j} e^{- \alpha (t-T_k^{(j)})} + \sum_{i\neq j} \sum_{ T_k^{(i)} < t} c_{il,j} e^{- \alpha (t-T_k^{(i)})}.$$
The initial conditions
$X_0,\lambda_0$ should be simulated according to the invariant distribution (and $\lambda_0^{(j)}$ should be larger than $\xi_j$). This measure of probability is not explicit. Thus we choose: $\lambda_0^{(j)}=\xi_j$ and $X_0=0$. Also, the exogenous intensities $\xi_j$ are chosen equal to $0.5$ for $j=1,2$. 
The weight matrix $c$ is chosen as:
 $$   c=\left(\begin{matrix}
    0.2 &0.3\\
    0.5 & 0.4\end{matrix}\right)$$
 and $\alpha=5$ (then the spectral radius of $H$ is nearly $0.02<1$).
    Figure \ref{fig:hawkes} shows one simulation of the Hawkes process for $T=10$. The graph on the left represents the jump times for the two components and the graph on the right represents the sum process of the intensities $\lambda_t^{(1)}+\lambda_t^{(2)}$.\\

Then we simulate $(X_\Delta, \dots X_{(n+1)\Delta})$ from an Euler scheme with $\delta= \Delta/5$. Indeed because of the additional term (when $a \neq 0$) to the best of our knowledge, it is not possible to use classical more sophisticated scheme. 
Finally we follow the following steps: 
\begin{enumerate}
\item Simulate $N^{(1)}, N^{(2)}$ on $[0,T]$ and deduce the jump times from the thinning method.
\item Create a large vector containing the grid of the $n$ times with time step $\Delta$ and the jump times, sorted. 
\item Simulate the process at these new times through an Euler scheme:\\
$X_{t_k}= (t_k-t_{k-1}) b( X_{t_k})+ \sqrt{t_k-t_{k-1}} \sigma(X_{t_k}) W+ a(X_{t_k})\one_{t_k = \text{'jump time'}} $.
\item Keep only in the grid the points corresponding to the regular grid of time step $\Delta$ without considering the jump times.
\end{enumerate}
A simulation algorithm is also detailed in \cite{DLL} Section 2.3.\\ 

We choose to show only relevant examples of processes with different types of drift function. There are four: 
\begin{enumerate}
\item Model 1: (Ornstein Uhlenbeck diffusion) $b(x)=-2x, ~\sigma(x)=1$, with $a(x)=\min(|x|, 5)$.
\item Model 2:  $b(x)=-(x-1/4)^3-(x+1/4)^3, ~\sigma(x)=1$, $a(x)=\min(|x|, 5)$.
\item Model 3:  $b(x)=-2x + \sin(3x), ~\sigma(x)=\sqrt{ \frac{3+x^2}{1+x^2}}$, $a(x)=\min(|x|, 5)$.
\item Model 4:  $b(x)=-2x, ~\sigma(x)=1$, $a(x)=0.2 x$.
\end{enumerate}

Model 1 is the simplest one. The second model do not satisfy the assumptions because in particular $b$ is not Lipschitz. The process is not ergodic and is really unstable. In order to obtain a result in this case we kept only the trajectories that do not explode in finite time. 
The third model has drift function similar to model 1, but with a larger diffusion coefficient. Finally, model 4 illustrate a case where the function $a$ is not positive and change sign. This last case can be more accurate to some situations (not neurons) where the impact of the network represented by the process $(N_t)$ is not only exiting.

Assumption \ref{ass:coeff}.2. is particularly strong and not verified in practice. Nevertheless it is true on the compact of estimation.

\begin{figure}
\centering
\includegraphics[height=8cm, width=11cm]{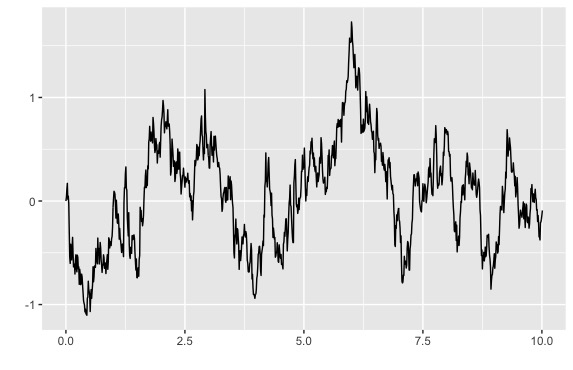}
\caption{Example of simulation of $X$ for model 1 with $n=100$ and $\Delta=1/10$, with the Hawkes process represented on Figure \ref{fig:estimb}.}
\label{fig:X}
\end{figure}

Figure \ref{fig:X} shows one trajectory of the process $(X_t)$ along time for model 1 when $n=100$ and $\Delta=1/10$, corresponding to the Hawkes process shown on Figure \ref{fig:hawkes}.

\subsection{Results of estimation of $b$}

\begin{figure}
\centering
\includegraphics[height=6cm, width=13cm]{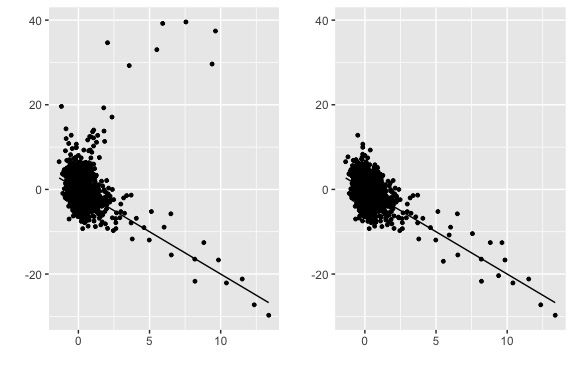}
\includegraphics[height=6cm, width=13cm]{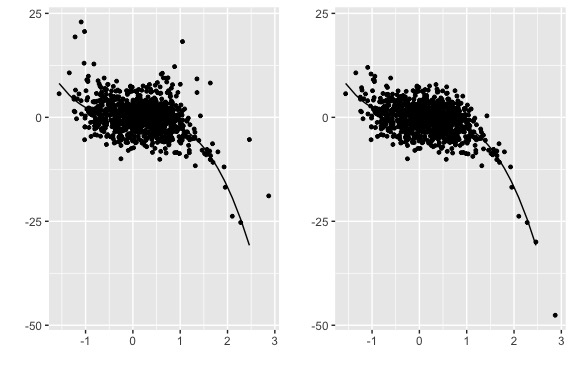}
\caption{At the left, $Y_{k\Delta} $ versus $X_{k\Delta}$, at the right, $U_{k\Delta}$ versus $X_{k\Delta}$. On the top, model 1, at the bottom, model 2. The true function $b$ is in plain line. Simulations for $n=1000$ and $\Delta=1/10$.}
\label{fig:XvsU}
\end{figure}

The results are illustrated for the trigonometric basis on the compact $A=[-1,1]$.
Then, the regression is done on the random variables $U_{k\Delta}= Y_{k\Delta}- T_{M, k \Delta}$. This requires to compute the term $T_{M,k\Delta}$ that we remind the reader is given by
$\displaystyle T_{M, k \Delta} = \frac{1}{\Delta} \int_{k\Delta}^{(k+1)\Delta}a(X_{s-}) \sum_{j=1}^Md{N}^{(j)}_s.$ This term is thus approximated by
$\displaystyle a(X_{k\Delta})\sum_{j=1}^M(N^{(j)}_{(k+1)\Delta}-N^{(j)}_{k\Delta})$. 
Figure \ref{fig:XvsU} illustrates the gain considering variables $U_{k\Delta}=Y_{k\Delta}-T_{M, k\Delta}$ instead of $Y_{k\Delta}$. For model $1$ and model $2$ (top and bottom) the left graph is $Y_{k\Delta}  $ versus $X_{k\Delta}$ together with the drift function $b$ in plain line, and on the right $U_{k\Delta}  $ versus $X_{k\Delta}$.
As attended, one can notice that the right graphs are more suited for regression, the points follows the curve of the true function $b$.

\begin{figure}
\centering
\includegraphics[height=8cm, width=10cm]{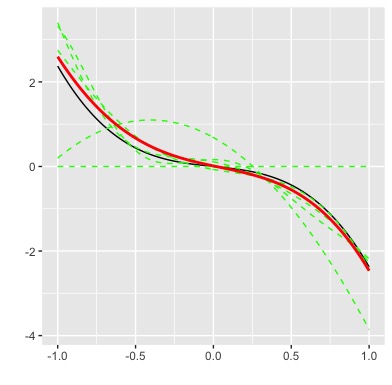}
\caption{Collection of estimators $\w{b}_m$ in (green) dotted line, true function $b$ in dark plain line and final estimator $\w{b}_{\w{m}}$ in red bold line for Model 2.}
\label{fig:estimb}
\end{figure}

Then, we have implemented the least squares estimator $\w{b}_m$ with the trigonometric basis. The dimension $D_m= 2m+1$ with $m \in \{1, \dots, m_{\text{max}}\}$ with $m_{\text{max}}= 20$ on simulations.
The adaptive procedure given in formula \eqref{eq:estimadapt} and \eqref{eq:pen} is implemented with the true value $\sigma_1$ (the bound of function $\sigma$ is known).
The calibration constant $\rho$ in the penalty function given in Equation \eqref{eq:pen} is calibrated on a large preliminary simulation study where 
we investigate various models with known parameters
and let $\rho$ vary. This constant
is chosen equal to $3$. 
Figure
\ref{fig:estimb} shows for model $2$ a collection of estimators $\w{b}_m$ in (green) dotted line, together with the true function $b$ in dark plain line and the chosen estimator in red bold line. 
Finally, we have computed a Monte-Carlo approximation of the empirical risks from $1000$ simulations. We show the results for $\Delta \in \{1/10,1/100\}$ with $n\in \{ 1000, 10000\}$.
On Table \ref{tab:emprisk} one can see the influence of the parameters $n$ and $\Delta$ on the empirical risk, and also the impact of the model. 
When $\Delta$ is multiplied by 10 and $n$ fixed the risk is also multiplied by $10$. But if $T=n\Delta$ is unchanged the influence of $\Delta$ is less clear. 
It is in line with the result of Theorem \ref{theo:adapt}.

The most difficult case is $T=10$. Of course in this case the process may not be stationary yet as we do not start from the invariant distribution. Model $1$ shows the best results. Indeed the model is simple and satisfies all the assumptions. Nevertheless the other results are not so far. 
One can compare for example with the results obtain for the estimation of the drift in the classical jump-diffusion in \cite{SCHMISSER3}, the results are very close.

\begin{table}
\centering
\begin{tabular}{r|rrrr}
& $\begin{matrix} n&=1000\\ \Delta&=1/10\end{matrix}$ & $\begin{matrix} n&=10000\\ \Delta&=1/10\end{matrix}$ &  $\begin{matrix} n&=1000\\ \Delta&=1/100\end{matrix}$ &  $\begin{matrix} n&=10000\\ \Delta&=1/100\end{matrix}$ \\
\hline
Model 1 & $0.068$ & $0.011$  & $1.917$ & $0.061$\\
Model 2 & $ 0.080$ & $0.007$ & $2.948$ & $0.104$\\
Model 3 & $ 0.161$ & $0.014$ & $4.528$ & $0.154$\\
Model 4 & $ 0.058$ & $0.005$ & $4.398 $ & $0.065$
\end{tabular}
\caption{Empirical risk $\E[ \|\w{b}_{\w{m}}-b\|_n^2]$ from 1000 repetitions of Monte-Carlo on $[-1,1]$.}
\label{tab:emprisk}
\end{table}

\section{Discussion}\label{sec:discussion}

In this paper, we provide a nonparametric estimator of the drift function in a diffusion process with jumps driven by a Hawkes process.
Section~\ref{sec:estimation} is devoted to theoretical guarantees for the non-adaptive and adaptive estimator.
Section~\ref{sec:simu} deals with the implementation of the process and of the estimator. 

In future works, we plan to generalize the obtained results. There are many remaining questions of interest. The case of an unknown diffusion coefficient $\sigma$ and an unknown jump coefficient $a$ are obviously crucial issues.
Indeed, the presented results require the knowledge of $\sigma_1^2$ (the square of the bound on volatility, see Assumption \ref{ass:coeff}) and of the coefficient $a(\cdot)$ (for the computation of $T_M$ in particular).
The relation 
$$ \frac{(X_{(k+1)\Delta}-X_{k\Delta})^2}{\Delta}=\sigma^2(X_{k\Delta})+a^2(X_{k\Delta}){\left( \frac{1}{\Delta} \int_{k\Delta}^{(k+1)\Delta} \sum_{j=1}^M{\lambda}^{(j)}_s ds \right)}+
\text{terms}$$
should lead to some estimation of $\sigma^2, a, \sigma^2+a $ (or at least in the parametric case). Nevertheless it is not a simple issue because the additional terms here are not centered and are hard to control. It will be investigated in future works.

One may also think of generalizing the initial model. For instance, the case of unknown kernel functions $h$ for the Hawkes process is of interest. Indeed, as several nonparametric estimators of $\lambda$ are available \citep[mainly][]{BACRYkernel,Hansen15} an estimator of $\lambda$ could be plug-in the procedure of the estimation of the coefficients.

Finally, this estimation procedure will be tested on a neuronal dataset and put in competition with classical methods of machine learning. This is a work in progress.

\section{Proofs}\label{sec:proofs}
\label{sec:proof}

We drive the reader attention on the fact that through the proofs 
the constants are generically noted $C,C'...$ 
for a strictly positive, real-valued constant, not necessarily the same.
Sometimes, the notation $C(k)$ is used to emphasize that the constant depends
on $k$ in some unspecified manner.

\subsection{Preliminary result}
Let us remind the reader a result obtained in \cite{BDH}, named Lemma 4 and Lemma 6. 
\begin{lemme}\label{lem:suitebacry}
For all $\delta >0$, 
\begin{enumerate}[label=(\roman*)]
\item there exists a positive constant $C$ such that
\begin{equation}\label{eq:bacry2}
 \sum_{j=1}^M\E[(N_{t+\delta}^{(j)}-N_t^{(j)})^2] \leq C (\delta+ \delta^2),
\end{equation}

\item
and there exists a positive constant $C$ such that
\begin{equation}\label{eq:bacrylem2}
 \sum_{j=1}^M\E[(N_{t+\delta}^{(j)}-N_t^{(j)})^4] \leq C (\delta+ \delta^2+\delta^4).
\end{equation}
\end{enumerate}
\end{lemme}


%
%
\n \textit{Proof of Lemma \ref{lem:suitebacry}}. 
 \textit{(i)} Following \cite{BACRY}, 
denoting $A^{(j)}_t:= N^{(j)}_t-\E[N^{(j)}_t]$, it comes
\begin{eqnarray*}
 \sum_{j=1}^M\E[(N_{t+\delta}^{(j)}-N_t^{(j)})^2] &\leq &2 \sum_{j=1}^M\E[(A_{t+\delta}^{(j)}-A_t^{(j)})^2] + (\E[(N_{t+\delta}^{(j)}-N_t^{(j)} ])^2\\
&\leq & 2 \|\E[ N_{t+\delta}-N_t]\|^2 + 2 \E[ \|A_{t+h}-A_t\|^2]
\end{eqnarray*}
where $\|.\|$ is the Euclidean norm on $\mathbb{R}^M$. The first term is controlled by Equation \eqref{eq:bacry}: $\E[N_{t+\delta}-N_t] = \delta (I_d-H)^{-1}\xi$ where $\xi$ is the M-vector of $\xi_j$'s and $H=( c_{j,k}/\alpha )_{j,k \in \{1, \dots, M\}}$, it is of order $\delta^2$. The second term is controlled by Lemma 6 of \cite{BDH} and is of order $\delta$. There are constants that depend on $c_{j,k}, \alpha, \xi_i$. \\


 \textit{(ii)} The same notations as in \cite{BACRY} are used:
\[
A_t^{(j)}=N_t^{(j)}-\E[N_t^{(j)}]=\widetilde N^{(j)}_t+\int_0^t\sum_{i=1}^M\psi_{ij}\widetilde N_s^{(j)}ds, 
\]
where $\widetilde N_t^{(j)}=N_t^{(j)}-\int_0^t\lambda_t^{(j)}$ is a local martingale, $\psi_{ij}=\sum_{n\geq 1}h_n^{ij}$ and $h_n^{ij}$ is such that $h_{n+1}^{ij}(t)=\int_0^th_{ij}(t-s)h_n^{ij}(s)ds$ and $h_1^{ij}=h_{ij}$. With these notations, we have for a constant $C$ that may change from a line to another
\[
\E\left[(N_{t+\delta}^{(j)}-N_t^{(j)})^4\right]\leq C\E \left[(A_{t+\delta}^{(j)}-A_t^{(j)})^4\right]+C\E\left[N_{t+\delta}^{(j)}-N_t^{(j)}\right]^4.
\]
From Equation \eqref{eq:bacry}, we have $\E[N_{t+\delta}^{(j)}-N_t^{(j)}]^4\leq C\delta^4$. The second term in the previous inequality can be controlled as follows
\begin{equation*}
\E[(A_{t+\delta}^{(j)}-A_t^{(j)})^4]\leq C\E[(\widetilde{N}_{t+\delta}^{(j)}-\widetilde N_{t}^{(j)})^4]+C\E\left[\left(\int_0^{+\infty}\sum_{i=1}^M\psi_{ij}(u)(\widetilde{N}_{t+\delta-u}^{(j)}-\widetilde N_{t-u}^{(j)})du\right)^4\right].
\end{equation*}
From \cite{BACRY}, Lemma 9, we have 
\[
\E\left[\underset{t\leq s\leq t+\delta}{\sup}(\widetilde N_s^{(j)}-\widetilde N_t^{(j)})^4|\mathcal{F}_t\right] \leq C(\delta+\delta^2).
\]
We deduce that 
\begin{align*}
\E[(A_{t+\delta}^{(j)}-A_t^{(j)})^4]&\leq C(\delta+\delta^2)+C\int_0^{+\infty}\Big(\sum_{i=1}^M\psi_{ij}(u)\Big)^4\E[(\widetilde{N}_{t+\delta-u}^{(j)}-\widetilde{N}_{t-u}^{(j)})^4]du\\
&\leq C(\delta+\delta^2).
\end{align*}
The last inequality holds due the fact that $\int_0^\infty\Big(\sum_{i=1}^M\psi_{ij}(u)\Big)^4du<+\infty$ for every $j=1,\dots,M$.
Indeed, for $h_{i,j}=c_{i,j}e^{-\alpha t}$ for $i,j\in \{1,\dots,M\}$, we have
\[
h_n^{i,j}(t)=c_{i,j}^n\dfrac{t^{n-1}}{(n-1)!}e^{-\alpha t}\quad \text{and}\quad  \psi_{ij}(t)=\sum_{n\geq 1}h_n^{ij}(t)=c_{i,j}e^{-(\alpha-c_{i,j})t}.
\]
Thus, if $\alpha\geq c_{i,j}$ for all $i,j\in\{1,\dots,M\}$, we have $\int_0^\infty\Big(\sum_{i=1}^M\psi_{ij}(u)\Big)^4du<+\infty$ for every $j=1,\dots,M$. $\Box$

\subsection{Proof of Proposition \ref{resultatcle}}

Following \cite{GLOTER} Proposition 5.1,
we first show that for $p=2$ or $p=4$, $\displaystyle \E \left[
\underset{s\in  [t, t+\delta]}{\sup} |X_s|^{p} \Big\vert {\mathcal{F}}_t
\right]\leq c(p)(1+|X_t|^p)$. Choose $s\in [t, t+\delta]$ with $\delta <1$.\\

\n Then,
\begin{eqnarray*}
X_s-X_t= \int_t^s b(X_u)du+\int_t^s \sigma(X_u)dW_u+ \sum_{j=1}^M \int_t^s a(X_{u-}) dN^{(j)}_u.
\end{eqnarray*}

Using Burkholder–Davis–Gundy inequality, the fact that $\sigma$ and $b$ are Lipschitz and that $a$ is bounded, we have for $p=2$ or $p=4$
\begin{align*}
\displaystyle \E \left[\underset{s\in  [t, t+\delta]}{\sup} |X_s|^{p} \Big\vert {\mathcal{F}}_t\right]&\leq c(p)|X_t|^p+c(p)\E\Big[\Big(\int_t^s\sigma^2(X_u)du\Big)^{p/2}\Big|\mathcal{F}_t\Big]\\
&+c(p)\E\Big[\Big(\int_t^s|b(X_u)|du\Big)^p\Big|\mathcal{F}_t\Big]+c(p)\E\Big[\Big(\sum_{j=1}^M\int_t^sa(X_{u-})dN_u^{(j)}\Big)^p\Big|\mathcal{F}_t\Big]\\
&\leq c(p)|X_t|^p+c(p)\E\Big[\Big(\int_t^s(|\sigma(X_u)|^p+|b(X_u)|^p)du\Big|\mathcal{F}_t\Big]\\
&+c(p,M)a_1^p\sum_{j=1}^M\E[(N_s^{(j)}-N_t^{(j)})^p|\mathcal{F}_t]\\
&\leq c(p)|X_t|^p+c(p)\int_t^s\E[(1+|X_u|^p)|\mathcal{F}_t]du+c(p,M)a_1^p\sum_{j=1}^M\E[(N_s^{(j)}-N_t^{(j)})^p|\mathcal{F}_t],
\end{align*}
where the constants may change from a line to another. 
We then introduce the following notation: $\phi(s)\leq \underset{u\in[t,s]}{\sup}\E[|X_u|^p|\mathcal{F}_t]$.
We have with this notation
\[
\phi(s)\leq c(p,M)(1+|X_t|^p+a_1^p\sum_{j=1}^M\E[(N_s^{(j)}-N_t^{(j)})^p|\mathcal{F}_t])+c(p)\int_t^s\phi(u)du. 
\]
From Assumption \ref{ass:coeff}.1, $\phi(s)$ is almost surely finite. We can apply Gronwall's Lemma (see in Appendix Lemma \ref{th:Gronwall}) to obtain
\begin{align*}
\phi(s)&\leq c(p,M)\Big(1+|X_t|^p+a_1^p\sum_{j=1}^M\E[(N_s^{(j)}-N_t^{(j)})^p|\mathcal{F}_t]\Big)+c(p)(1+|X_t|^p)(s-t)\\
&+c(p,M)a_1^p\int_t^s\sum_{j=1}^M\E[(N_u^{(j)}-N_t^{(j)})^p|\mathcal{F}_t]e^{c(p)(t-u)}du.
\end{align*}
Since $s\in[t,t+\delta]$ and $\delta<1$, applying Lemma \ref{lem:suitebacry}, we obtain 
\begin{equation}
\label{eq:gronwall}
\underset{u\in[t,t+\delta]}{\sup}\E[|X_u|^p|\mathcal{F}_t]\leq c(p)(1+|X_t|^p).
\end{equation}

\n Now, we consider 
$\Delta_{t, t+\delta}:= \underset{s \in [t, t+\delta]}{\sup}|X_s-X_t|$. As previously, we have
\begin{eqnarray*}\label{eq:eqresultatc}
\E\left[(\Delta_{t, t+\delta})^p|{\mathcal{F}}_t \right]
&\leq &  c(p) 
\E \left[ 
\left( 
\int_t^{t+\delta} \sigma^2(X_u)du 
\right)^{p/2} \Big \vert
{\mathcal{F}}_t
\right]+ 2 \E \left[ 
\left( 
\int_t^{t+\delta} |b(X_u)|du 
\right)^{p} \Big \vert
{\mathcal{F}}_t
\right]
 \nonumber \\
 &&+  c(p,M,a_1)\sum_{j=1}^M
 \E \left[ 
\left( N^{(j)}_{t+ \delta}-N^{(j)}_{t} \right)^{p} \Big \vert\bar{\mathcal{F}}_t
\right]\\
&\leq & c(p)
\E \left[ 
\delta^{p/2}
\sup_{t\leq s \leq  t+\delta} \sigma^p(X_s) 
+\delta^p
\sup_{t\leq s \leq  t+\delta} |b(X_s)|^p 
 \Big \vert
 {\mathcal{F}}_t
\right]
 \nonumber \\
 &&+  c(p,M,a_1)\sum_{j=1}^M\E \left[ 
\left( N^{(j)}_{t+ \delta}-N^{(j)}_{t} \right)^{p} \Big \vert \mathcal{F}_t
\right]
\end{eqnarray*}
and
\begin{eqnarray*}
\E\left[(\Delta_{t, t+\delta})^p|{\mathcal{F}}_t \right]
& \leq &  c(p)\delta^{p/2} ~\E \left[\sup_{t\leq s \leq  t+\delta}(1+|X_s|^p)\Big \vert
{\mathcal{F}}_t\right]+ c(p,M,a_1)\sum_{j=1}^M
\E\left[\left( N^{(j)}_{t+ \delta}-N^{(j)}_{t} \right)^{p} 
\Big \vert
{\mathcal{F}}_t\right] .
\end{eqnarray*}

Finally, according to Equations \eqref{eq:bacry2}, \eqref{eq:bacrylem2} and \eqref{eq:gronwall}, since $p=2$ or $p=4$, there exists a constant $C(M,a_1)$, such that 
\begin{eqnarray*}
\E\left[(\Delta_{t, t+\delta})^p\right]
& \leq &C(M,a_1)\delta(1+\E[|X_t|^p]).
\end{eqnarray*}

Finally 
since
$\E[X_t^p]< \infty$, we deduce the assertion.
$\Box$

\subsection{Additional results}

\begin{lemme}
\label{lem:tec}
For $p=2$ or $p=4$, it yields
\begin{align*}
\mathbb{E}\Big[\Big(\int_{k\Delta}^{(k+1)\Delta}b(X_{k\Delta})-b(X_s)ds\Big)^p\Big]
\leq C(M, a_1)\Delta^{3p/2}.
\end{align*}
\end{lemme}
\n \textbf{Proof of Lemma \ref{lem:tec}} 
According to Cauchy-Schwarz and Jensen's inequalities, we have
\begin{eqnarray*}
\left( \int_{k\Delta}^{(k+1)\Delta} b(X_{k\Delta})-b(X_s)ds \right)^p
&=& \Delta^p\left( \frac{1}{\Delta}\int_{k\Delta}^{(k+1)\Delta} b(X_{k\Delta})-b(X_s)ds \right)^p\\
& \leq & \Delta^p \left( \frac{1}{\Delta} \int_{k\Delta}^{(k+1)\Delta} ( b(X_{k\Delta})-b(X_s))^p ds \right)\\
& \leq &  \Delta^{p-1}  \int_{k\Delta}^{(k+1)\Delta} ( b(X_{k\Delta})-b(X_s))^p ds.
\end{eqnarray*}
Applying Proposition \ref{resultatcle},
gives
\begin{align*}\E\left[ \left( \int_{k\Delta}^{(k+1)\Delta} b(X_{k\Delta})-b(X_s)ds \right)^p\right]&\leq 
\Delta^{p-1} \E \left[   \int_{k\Delta}^{(k+1)\Delta} (b(X_{k\Delta})-b(X_s))^pds    \right]\nonumber\\
&\leq C(M, a_1) \Delta^{p+1}. \quad \Box
\end{align*}

\begin{lemme}
\label{lem:esp}
For some constant $C(M, a_1)>0$ (depending on $M$ and $a_1$) we have
$$\E[ I_{k\Delta}^2] \leq C(M, a_1)\Delta, \quad \E[I_{k\Delta}^4]\leq C(M, a_1)\Delta.$$
There exists $C>0$ such that $$\E[Z_{k\Delta}^4\vert{\mathcal{F}}_{k\Delta}] \leq C\dfrac{\sigma_1^4}{\Delta^2}.$$
\end{lemme}

\n \textbf{Proof of Lemma \ref{lem:esp}}\\
The first point comes from Lemma 
 \ref{lem:tec} with $p=2$,
the second point can be deduced when $p=4$.
\n Finally, with Burkholder–Davis–
Gundy inequality we find,
$$\E \left[Z_{k\Delta}^4 \vert {\mathcal{F}}_{k\Delta} \right]
\leq \frac{C}{\Delta^4} \E \left[\left( \int_{k\Delta}^{(k+1)\Delta} \sigma^2(X_s) ds\right)^2\Big|{\mathcal{F}}_{k\Delta} \right] \leq \frac{C}{\Delta^4}\sigma_1^4 \Delta^2 = C \frac{\sigma_1^4}{\Delta^2}.\quad \Box
$$

\subsection{Proof of Proposition \ref{prop:MISE}}
From the definition of the contrast function it comes for $t\in \mathbb{L}^2(A)$,
\begin{eqnarray*}
\gamma_{n,M}(t)-\gamma_{n,M}(b) &=& \| t-b\|^2_n +\frac{2}{n} 
 \sum_{k=1}^{n}  (b-t)(X_{k\Delta}) I_{k\Delta} + \frac{2}{n} 
 \sum_{k=1}^{n}  (b-t)(X_{k\Delta}) Z_{k\Delta}.
\end{eqnarray*}
Let us define
%
\begin{equation}\label{mu}
\mu_n(t):=  \frac{1}{n} \sum_{k=1}^{n}  t(X_{k\Delta})Z_{k\Delta}.
\end{equation}
By definition of $\w{b}_m$, we have
$$\gamma_{n,M}(\w{b}_m)-\gamma_{n,M}(b)\leq \gamma_{n,M}({b}_m)-\gamma_{n,M}(b),
$$
where $b_m \in \mathcal{S}_m$, for example $b_m$ is the orthogonal projection of $b$ on $\mathcal{S}_m$.
This implies 
\begin{eqnarray*}
\| \w{b}_m-b\|_n^2 &\leq & \| {b}_m-b\|_n^2+ 
 {\frac{2}{n} \sum_{k=1}^{n}  ( \w{b}_m-b_m)(X_{k\Delta}) I_{k\Delta} }+  2 \mu_n( \w{b}_m-b_m)
\end{eqnarray*}
(where $b$ denotes the restriction of $b$ to $A$).
Let us denote:
$$ \mathcal{B}_{m,\pi}= \{ t\in \mathcal{S}_m, \|t\|_{\pi^X} =1\},$$
 If $t$ is a deterministic function,  $\E[\|t\|^2_n]=\|t\|_{\pi^X}^2$. \\

\n Using Cauchy-Schwarz inequality and the relation 
$2xy \leq x^2/d + dy^2$ for $x, y, d>0$, it comes
\begin{eqnarray*}
\| \w{b}_m-b\|_n^2 &\leq & \| {b}_m-b\|_n^2+  2\| \w{b}_m-b_m\|_{n} \left( \frac{1}{n} \sum_{k=1}^{n} I_{k\Delta} ^2 \right)^{1/2}+ 2\| \w{b}_m-b_m\|_{\pi^X}
\underset{t\in \mathcal{B}_{m, \pi^X}}{\sup} \vert\mu_n(t)|\\
& \leq &  \| {b}_m-b\|_n^2 +\frac{1}{d}\| \w{b}_m-b_m\|_{\pi^X}^2+
d \underset{t\in \mathcal{B}_{m, \pi^X}}{\sup} \vert\mu_n(t) \vert ^2 +  \frac{1}{d} \| \w{b}_m-b_m\|_{n}^2+ \frac{d}{n} \sum_{k=1}^{n} I_{k\Delta} ^2.
\end{eqnarray*}

%
\n First, let us consider the following set where the $\pi^X$-norm and the empiric norm are equivalent:
\begin{equation}
\Omega_n:= \left\{\omega, \forall t \in \mathcal{S} \backslash \{0\},~ \left\vert \frac{\|t\|^2_{n}}{\|t\|^2_{\pi^X}}-1  \right\vert \leq 1/2 \right\}.
\end{equation}
\n \textbf{Bound of the risk on $\Omega_n$}\\
On $\Omega_n$, we have $\|\w{b}_m - b_m \|_{\pi^X}^2 \leq 2 \|\w{b}_m-b_m\|_n^2$ and $\|\w{b}_m-b_m\|_n^2 \leq 2( \|\w{b}_m-b\|_n^2+\|{b}_m-b\|_n^2)$ 
thus 
 it comes 
\begin{equation}\label{eq:boudonomegan}
\| \w{b}_m-b\|_n^2 \leq  13\| {b}_m-b\|_n^2
+49\frac{1}{n} \sum_{k=1}^{n} I_{k\Delta} ^2+ 49\underset{t\in \mathcal{B}_{m, \pi^X}}{\sup} \vert\mu_n(t) \vert ^2 
\end{equation}
Denote $(\psi_{\ell})_{\ell}$ an orthogonal basis of $\mathcal{S}_m$ for the $\mathbb{L}^2_{\pi^X}$-norm (thus $\displaystyle \int \psi^2_{\ell}(x)\pi^X(x)dx=1$). 
It comes from Cauchy-Schwarz inequality that,
\begin{eqnarray}\label{eq:decompbase}
\underset{t\in \mathcal{B}_{m, \pi^X}}{\sup} \vert\mu_n(t) \vert ^2= \underset{\sum_\ell \alpha_\ell^2 \leq 1}{\sup}\mu^2_{n,M}\left(\sum_\ell \alpha_\ell \psi_\ell \right)  \leq \underset{\sum_\ell \alpha_\ell^2 \leq 1}{\sup} \left(  \sum_{\ell}\alpha_\ell^2    \right) \left( \sum_{\ell}\mu_n^2(\psi_\ell) \right)=\sum_{\ell=1}^{D_m}\mu_n^2(\psi_\ell).
\end{eqnarray}
\n Moreover, 
\begin{eqnarray}\label{eq:boundmu}
\E \left[\underset{t\in \mathcal{B}_{m, \pi^X}}{\sup}\mu_n^2(t)  \right] &\leq & \sum_{ \ell=1}^{D_m}\E[\mu_n^2(\psi_\ell)]\nonumber\\
& \leq & \frac{2}{n^2 \Delta^2} \sum_{k=1}^{n} \E \left[    \sum_{ \ell=1}^{D_m } \psi_\ell^2(X_{k\Delta})
 \E\Big[\Big(\int_{k\Delta}^{(k+1)\Delta}\sigma(X_s)dW_s\Big)^2|{\mathcal{F}}_{k\Delta}\Big]   \right]\nonumber \\
&\leq & \frac{2}{n^2 \Delta^2} \sum_{k=1}^{n} \E \left[    \sum_{ \ell=1}^{D_m } \psi_\ell^2(X_{k\Delta}) \E\Big[\int_{k\Delta}^{(k+1)\Delta}\sigma^2(X_s)ds|{\mathcal{F}}_{k\Delta}\Big]   \right]\nonumber \\
& \leq & \frac{2D_m  \sigma_1^4}{n\Delta} 
\end{eqnarray}
(if $\sigma$ is not bounded one can use the stationarity of the process). 
With Lemma \ref{lem:esp} the obtained result is
\begin{equation}\label{eq:boundonomega2}
\E[\| \w{b}_m-b\|_n ^2\mathds{1}_{\Omega_n}] \leq C_1 \E[\|b-b_m \|_n^2]+C_2\frac{D_m}{n\Delta}+C_3\Delta
\end{equation}
with $C_1= 13$, $C_2= 98\sigma_1^4$, $C_3$ depending on $a_1,M,\sigma_1$. 
\\

\n \textbf{Bound of the risk on $\Omega_n^c$ }\\

Let us set $e=(e_\Delta,\dots,e_{n\Delta})$, where $e_{k\Delta}:=U_{k\Delta}-b(X_{k\Delta})=I_{k\Delta}+Z_{k\Delta}$ and $\Pi_mU=\Pi_m(U_\Delta,\dots,U_{n\Delta})=(\w{b}_m(X_\Delta),\dots,\w{b}_m(X_{n\Delta}))$ where $\Pi_m$ is the Euclidean orthogonal projection over $\mathcal{S}_m$. Then according to the projection definition,
\begin{align*}
\|\w{b}_m-b\|_n^2&=\|\Pi_mU-b\|_n^2=\|\Pi_mb-b\|_n^2+\|\Pi_m U-\Pi_mb\|_n^2\\
&\leq \|b\|_n^2+\|U-b\|_n^2=\|b\|_n^2+\|e\|_n^2.
\end{align*}
Cauchy-Schwarz inequality implies,
\begin{align*}
\E[\|e\|^4_n]&=\E\Big[\Big(\dfrac{1}{n} \sum_{k=1}^{n}(I_{k\Delta}+Z_{k\Delta})^2\Big)^2\Big]\\
&\leq \dfrac{1}{n} \sum_{k=1}^{n}2^3\E\Big[I_{k\Delta}^4+Z_{k\Delta}^4\Big].
\end{align*}
With the controls given in Lemma \ref{lem:esp}, the additional following result will end the proof.
\begin{lemme}\label{lem:omegacproba}
If $n\Delta / \ln^2(n) \rightarrow \infty $ if $D_n \leq O((n\Delta)/(\ln^2(n)))$, the probability of the event $\Omega_n^c$ is bounded as follows:
$$\mathbb{P}(\Omega_n^c) \leq \frac{c_0}{n^4}.$$
\end{lemme}
\n This Lemma is proven in Subsection \ref{sec:lemomegacproba}.
\\

\n From Lemma \ref{lem:esp}, Lemma \ref{lem:omegacproba}, Cauchy-Schwarz inequality, we finally obtain 
$$\E[ \| e\|_n^2 \one_{\Omega_n^c}] \leq \left( \E[ \| e\|_n^4] ^{1/2} \right) \left(\P(\Omega_n^c)^{1/2} \right) \leq  \dfrac{C}{n\Delta}.
$$
Using that $\|\w b_m-b\|_n^2\leq \|b\|_n^2+\|e\|_n^2$, we deduce under the stationarity assumption, that
\begin{eqnarray*}
\E[\|\hat b_m-b\|_n^2\mathds{1}_{\Omega_n^c}] &\leq &
\E \left[ \frac{1}{n} \sum_{k=1}^{n} b(X_{k\Delta})^2   \one_{\Omega_n^c}\right]
+ \E \left[ \frac{1}{n} \sum_{k=1}^{n} (e_{k\Delta})^2   \one_{\Omega_n^c}\right]\\
&\leq & \E \left[ b^2(X_0)  \one_{\Omega_n^c} \right] + 
\E \left[ e^2_{\Delta}  \one_{\Omega_n^c} \right]\\
&\leq & \E \left[  b^4(X_0)\right]^{1/2} \P(\Omega_n^c)^{1/2}
+
 \E \left[  e^4_{\Delta}\right]^{1/2} \P(\Omega_n^c)^{1/2}.
\end{eqnarray*}
From the Lipschitz property of $b$, we get 
$ \displaystyle  \E[b^4(X_0)] \leq C(1+ \E[X_0^4]) < \infty
$, and from the Lemma \ref{lem:esp} $\E[e^4_\Delta]^{1/2} \leq C'/\Delta$, thus
\begin{equation}\label{eq:riskomeganc}
E[\|\hat b_m-b\|_n^2\mathds{1}_{\Omega_n^c}]
\leq \frac{C''}{n\Delta}.
\end{equation}
Gathering Equation (\ref{eq:boundonomega2}) and (\ref{eq:riskomeganc}) gives the attended result. $\Box$\\

\subsection{Proof of Theorem \ref{theo:adapt}}

\n On $\Omega_n^c$, the proof can be lead as in the proof of Proposition \ref{prop:MISE}.\\

\n On $\Omega_n$, by definition of $\hat m$, we have $\gamma_n(\hat b_{\hat m})+\pen(\hat m)\leq \gamma_n(\w b_m)+\pen(m)$. We deduce with the same computation as previously leading to \eqref{eq:boundonomega2} that we have for all $m\in \mathcal{M}_n$,

\begin{eqnarray*}
\E[\| \w{b}_{\w m}-b\|_n^2 \mathds{1}_{\Omega_n}]
& \leq & 13  \| {b}_m-b\|_n^2 + 49 C \Delta + 49 \E\left[\underset{t\in \mathcal{B}_{m, \w m,\pi^X}}{\sup} \vert\mu_n(t) \vert ^2 \mathds{1}_{\Omega_n}\right]\\ 
&+&7\pen(m)-7\E[\pen(\w m)],
\end{eqnarray*}
where $\mathcal{B}_{m, m', \pi^X}=\{h\in S_m+S_{ m'}:\|h\|_{\pi^X}\leq 1\}$ and $\mu_n$ given in Equation \eqref{mu}. The challenge here is to compute the expectation of the supremum on a random ball. 
Let us introduce the notation
\begin{align*}
G_m(m')&:=\underset{t\in \mathcal{B}_{m, m', \pi^X}}{\sup}~|\mu_n(t)|.
\end{align*}
Then, 
\begin{eqnarray*}
\E[\| \w{b}_{\w m}-b\|_n^2 \mathds{1}_{\Omega_n}]
& \leq & 13 \E[\| {b}_m-b\|_n^2] + 49 C \Delta + 49 \E\left[G_m^2(\w{m})\mathds{1}_{\Omega_n}\right]\\ 
&+&7\pen(m)-7\E[\pen(\w m)]\\
& \leq & 13 \E[\| {b}_m-b\|_n^2] + 49 C \Delta + 49\sum_{m'\in \mathcal{M}_n}\E\left[G_m^2(m')-p(m,m'))\mathds{1}_{\Omega_n}\right]_+\\ 
&+& 49p(m, \w{m})+7\pen(m)-7\E[\pen(\w m)]
\end{eqnarray*}
where $\pen(m)$ must satisfy $7p(m,m')\leq \pen(m)+\pen(m')$.

To control the term $\displaystyle \E\left[G_m^2(m')-p(m,m'))\mathds{1}_{\Omega_n}\right]_+$  we have to prove a Bernstein-type inequality.
Here we follow \cite{BCV2}.
 The idea is presented in \cite{CGCR} Lemme 2 page 533. 
We have the following control. 
\begin{lemme}\label{lem:bernstein}
For any positive numbers $\varepsilon, ~v$, we have
$$\mathbb{P}\left( \sum_{k=1}^n t(X_{k\Delta})Z_{k\Delta} \geq n \varepsilon, ~  \|t\|_n^2 \leq v^2\right) \leq \exp\left( -\frac{n\varepsilon^2 \Delta}{2v^2 \sigma_1^2} \right).$$
\end{lemme}
The proof is based on the martingale $\int_0^s H_u dW_u, ~  H_u= \sum_{k=1}^n \one_{[k\Delta, (k+1)\Delta [}(u) t(X_{k\Delta})\sigma(X_u)$ (for the filtration $\sigma(X_u, ~ u \leq n \Delta)$), and follows exactly the one detailed in \cite{CGCR} and is omitted here.

Then, the chaining argument of \cite{BCV2} Proposition 6.1 page 45 can be applied here. We finally obtain:
$$
\E\left[G_m^2(m')-p(m,m'))\mathds{1}_{\Omega_n}\right]_+ \leq c \sigma_1^2 \frac{e^{-D_{m'}}}{n\Delta}
$$
and $p(m,m')= \rho \sigma^2_0 \frac{D_m+D_m'}{n\Delta}$, thus $\pen(m) \geq 7\rho \sigma_1^2 D_m / (n\Delta)$.
Finally, we have shown that
\begin{eqnarray*}
\E[\| \w{b}_{\w m}-b\|_n^2 \mathds{1}_{\Omega_n}]
& \leq & 13\E[ \| {b}_m-b\|^2_n] + 49 C \Delta + 49c \sigma_1^2 \sum_{m'\in \mathcal{M}_n} \frac{e^{-D_{m'}}}{n\Delta}+7\pen(m),
\end{eqnarray*}
 which ends the proof. $\Box$

\subsection{Proof of Lemma \ref{lem:omegacproba}}\label{sec:lemomegacproba}

The proof of this result follows the proof of Lemma 1 of \cite{CGCR}. 
The difficulty occurs here in the fact that the projection of the invariant measure $\pi$ onto the $X-$ coordinate $\pi^X$ is not bounded from below on the compact $A$. Then, the proof must be adapted from \cite{BCV}.

Under the exponential $\beta-$mixing control, similarly to \cite{CGCR} using Berbee's coupling method (for the random variables $X_{k\Delta}$), we obtain the decomposition :
\begin{equation}\label{eq:decP}
 \mathbb{P}(\Omega_n^c) \leq \mathbb{P}(\Omega_n^c \cap \Omega^*)+ \mathbb{P}(\Omega^{*c})
 \end{equation}
and $$\mathbb{P}(\Omega^{*c})\leq n \beta_X(q_n \Delta)$$
where $q_n$ is an integer such that $q_n < n $.
The notations are those used in \cite{CGCR}, in particular the notation $\text{}^*$ referring to the coupling variables is not detailed here. 

 Let us control the first term now.
Let $(\varphi_j)_j=1, \dots, D_n$ be an $\mathbb{L}^2(A)$-orthonormal basis of $\mathcal{S}_n$ and define the matrices:
$$V_{\pi^X}= \left[\left( \int_A \varphi_j^2(x) \varphi_{j'}^2(x) \pi^X(x)dx \right)^{1/2} \right]_{j,j'\in \{1, \dots, D_n\}^2}, \quad B= \left( \|\varphi_j \varphi_{j'}\|_\infty\right)_{j,j'  \in \{1, \dots, D_n\}^2}.$$
Besides, denote $L_n(\varphi):= \max( \rho^2(V_{\pi^X}), \rho(B))$, where for a symmetric matrix $A$:\\
$$\rho(A):= \sum_{a_j, \sum_{j}a_j^2 \leq 1} \sum_{j,j'} |a_j||a_j'| |A_{j,j'}|.$$
We note $\mathbb{P}^* = \mathbb{P}(\cdot \cap \Omega^*)$. Let us consider $\nu_n(t):=(1/n) \sum_{k=1}^n t(X_{k\Delta})- \E[t(X_{k\Delta})]$, $\mathcal{B}_{\pi^X}(0,1)= \{t \in \mathcal{S}_n, ~\|t\|_{\pi^X}\leq 1\}$ and $\mathcal{B}(0,1)= \{t \in \mathcal{S}_n, ~\|t\| \leq 1 \}$, where $\|.\|$ is the $\mathbb{L}^2$-norm. 
As on $A$, $\pi_0 \leq \pi^X(x)$, we have: 

$$\sup_{t \in \mathcal{B}_{\pi^X}(0,1)} |\nu_n(t^2)| = \sup_{t \in \mathcal{S}_n \backslash \{0\}} \left| \frac{\|t\|_n^2 }{\|t\|_{\pi^X}^2} -1\right| \leq \pi_0^{-1} \sup_{t \in \mathcal{B}_{(0,1)}} |\nu(t^2)| .$$
Thus, 
\begin{eqnarray*}
\mathbb{P}^* \left( \sup_{t \in \mathcal{B}_{\pi^X}(0,1)} |\nu_n(t^2)|  \geq \rho_0 \right) &\leq &\mathbb{P}^* \left( \sup_{t \in \mathcal{B}(0,1)} |\nu_n(t^2)| \geq \pi_0 \rho_0 \right) \\
& \leq & \mathbb{P}^* \left( \sup_{ \sum_{j=1}^{D_n}a_j^2 \leq 1} \sum_{j,j'} |a_j a_j'| | \nu(\varphi_j \varphi_{j'})| \geq \pi_0 \rho_0 \right).
\end{eqnarray*}
On the set $\{ \forall (j,j')\in \{1, \dots, D_n\}^2, ~|\nu_n(\varphi_j \varphi_{j'})| \leq 2 V_{\pi^X, j, j'}(2x)^{1/2}+3 B_{j,j'}x\}$ we have 
$$\sup_{ \sum_{j=1}^{D_n}a_j^2 \leq 1} \sum_{j,j'} |a_j a_j'| | \nu_n(\varphi_j \varphi_{j'})|  \leq 2 \rho(V_{\pi^X}) (2x)^{1/2}+3 \rho(B)x.$$
Choosing $x= (\rho_0 \pi_0)^2/(2 L_n(\varphi))$ with $\rho_0= 4(1+ \sqrt{1+3\pi_0})/(3\pi_0)$ it comes
$$\sup_{ \sum_{j=1}^{D_n}a_j^2 \leq 1} \sum_{j,j'} |a_j a_j'| | \nu_n(\varphi_j \varphi_{j'})|  \leq \pi_0/2.$$
Thus we have obtained: 
\begin{eqnarray*}
\mathbb{P}^*(\Omega_n^c) &=& \mathbb{P}^* \left( \sup_{t \in \mathcal{B}_{\pi^X}(0,1)} |\nu_n(t^2)|  \geq 1/2 \right) \\
& \leq & \mathbb{P}^*\left( \{\forall(j,j') \in \{1, \dots, D_n\}^2, ~ | \nu_n(\varphi_j \varphi_{j'})|  \leq 2 V_{\pi^X, j, j'}(2x)^{1/2}+3 B_{j,j'}x\} \right).
\end{eqnarray*}
 Then, we have to bound this last probability term. This is done in \cite{BCV} Claim 6 of Proposition 7 using a Bernstein's inequality. Let us denote as in \cite{BCV} $Z^*_{l,k}$ the random variables to which the inequality is applied. As the projection of the invariant measure $\pi$ onto the $X-$ coordinate $\pi^X$ is not bounded from below,
we have only
 $$E[(Z^*_{l,k})^2] \leq \int_A \varphi_j^2(x) \varphi_{j'}^2(x) \pi^X(x)dx = V_{\pi^X, j, j'}$$
 and still $\|Z_{l,k}^*\|_\infty \leq B_{j,j'}$.
Then, it comes, 
$$ \mathbb{P}^*\left( \{\forall(j,j') \in \{1, \dots, D_n\}^2, ~ | \nu_n(\varphi_j \varphi_{j'})|  \leq 2 V_{\pi^X, j, j'}(2x)^{1/2}+3 B_{j,j'}x\} \right) \leq 2 D_n^2 \exp\left(-\frac{nx}{q_n} \right).$$
The chosen value for $x$ implies that there is a constant $C>0$ such that 
$$ \mathbb{P}^*\left( \{\forall(j,j') \in \{1, \dots, D_n\}^2, ~ | \nu_n(\varphi_j \varphi_{j'})|  \leq 2 V_{\pi^X, j, j'}(2x)^{1/2}+3 B_{j,j'}x\} \right) \leq 2 D_n^2 \exp\left(-C\frac{n}{q_n L_n(\varphi)} \right).$$
It remains to control $L_n(\varphi)$. Following Lemma 2 of \cite{BCV} we have for the trigonometric basis according to Cauchy-Schwarz:
\begin{eqnarray*}
\rho^2(V_{\pi^X}) &\leq & \sum_{j,j' \in \{1, \dots, D_n\}^2} \int_A \varphi_j^2(x) \varphi_{j'}^2(x)\pi^X(x) dx \\
& \leq &  \left\| \sum_{j=1}^{D_n} \varphi_j^2 \right\|_\infty \sum_{j'=1}^{D_n} \int_A \varphi_{j'}^2(x) \pi^X(x) dx \\
& \leq &\phi_{\pi^X}  \phi_1^2 D_n^2 
\end{eqnarray*}
under Assumption \ref{ass:basis}.1  and \ref{ass:basis}.3. $\rho(B) \leq \phi_1^2D_n$.
Thus, $L_n(\varphi) \leq \max(1,\phi_{\pi^X}) \phi_1^2 D_n^2$.
Now, following Equation \eqref{eq:betamix}, $\beta_X(q_n \Delta) \leq e^{-\theta q_n \Delta}$, decomposition \eqref{eq:decP} leads to
$$\mathbb{P}(\Omega_n^c) \leq 1/n^4 + 2n^2 \exp\left(-C' \frac{n\Delta}{\ln(n)D_n}\right),$$
with $q_n=\lfloor 5\ln(n)/ (\theta\Delta)\rfloor +1  $.
Finally, for the trigonometric basis, if $\displaystyle D_n^2 \leq \frac{C'n\Delta}{6\ln^2(n)}$ thus $D_n \leq C' \sqrt{n\Delta}/\ln(n)$ and 
 $\mathbb{P}(\Omega_n^c) \leq c/n^4$, if $n\Delta / \ln^2(n) \rightarrow \infty$.
$\Box$


\section{Appendix: theoretical results}\label{sec:appendix}

\subsection{Gronwall's Lemma}
\begin{lemme}\label{th:Gronwall}
Let $\phi$, $\psi$ and $y$ be three non negative continuous functions on $[a,b]$, that satisfy
\[
\forall t\in[a,b], \quad y(t)\leq \phi(t)+\int_a^t\psi(s)y(s)ds.
\]
Then
\[
\forall t\in[a,b], \quad y(t)\leq\phi(t)+\int_a^t\phi(s)\psi(s)\exp\Big(\int_s^t\psi(u) du\Big)ds.
\]
\end{lemme}

\subsection{Talagrand's inequality}
The following result follows from the Talagrand concentration inequality \citep[see][]{klein}.

\begin{theo}\label{Talagrand}
Consider $n \in \mathbb{N}^*$, $\mathcal{F}$ a class at most countable of measurable functions, and  $(X_i)_{i\in\{1,...,n\}}$ a family of real independent random variables.
One defines, for all $f\in\mathcal{F}$, 
$$\nu_n(f)  =\frac{1}{n} \sum_{i=1}^{n} (f(X_i)-\mathbb{E}[f(X_i)]).$$
Supposing there are three positive constants $M$, $H$ and $v$ such that 
$\underset{f\in\mathcal{F}}{\sup} \|f\|_{\infty} \leq M$,\\
 $\mathbb{E}[\underset{f\in\mathcal{F}}{\sup} |\nu_n(f)| ] \leq H$, and $\underset{f\in\mathcal{F}}{\sup} ({1}/{n}) \sum_{i=1}^{n}  \mathrm{Var}(f(X_i)) \leq v$, then for all $\alpha>0$, 
\begin{eqnarray*}
\mathbb{E}\left[ \left(  \underset{f\in\mathcal{F}}{\sup} |\nu_n(f)|^2-2(1+2\alpha)H^2 \right)_+ \right]& \leq & \frac{4}{b} \left( \frac{v}{n} \exp \left(-b \alpha \frac{n H^2}{v} \right) \right.\\
&+& \left. \frac{49M^2}{b C^2(\alpha)n^2} \exp\left(-\frac{\sqrt{2}b C(\alpha)\sqrt{\alpha}}{7}\frac{nH}{M} \right) \right )
\end{eqnarray*}
with $C(\alpha)=(\sqrt{1+\alpha}-1) \wedge 1$, and $b=\frac{1}{6}$.
\end{theo}

\section*{Acknowledgments}

The authors are very greatfull to CNRS for the financial support PEPS.
They also thank Gilles Pages for the fruitful discussions.
The authors wish to thank Patricia Reynaud Bouret for her advice and her listening. 
Finally, the authors are very grateful to Eva L\"ocherbach who supported the project and helped to improve the paper.

\newpage

\bibliographystyle{ScandJStat}
\bibliography{BIB}

\begin{thebibliography}{58}
\providecommand{\natexlab}[1]{#1}
\providecommand{\url}[1]{\texttt{#1}}
\providecommand{\urlprefix}{URL }
\expandafter\ifx\csname urlstyle\endcsname\relax
  \providecommand{\doi}[1]{doi:\discretionary{}{}{}#1}\else
  \providecommand{\doi}{doi:\discretionary{}{}{}\begingroup
  \urlstyle{rm}\Url}\fi

\bibitem[{Amorino \& Gloter(2018)}]{amorino2018contrast}
Amorino, C. \& Gloter, A. (2018).
\newblock Contrast function estimation for the drift parameter of ergodic jump
  diffusion process.
\newblock \emph{arXiv preprint arXiv:1807.08965} .

\bibitem[{Bacry \emph{et~al.}(2012)Bacry, Dayri \& Muzy}]{BACRYkernel}
Bacry, E., Dayri, K. \& Muzy, J. (2012).
\newblock Non-parametric kernel estimation for symmetric hawkes processes.
  application to high frequency financial data.
\newblock \emph{THE EUROPEAN} .

\bibitem[{Bacry \emph{et~al.}(2013)Bacry, Delattre, Hoffmann \& Muzy}]{BDH}
Bacry, E., Delattre, S., Hoffmann, M. \& Muzy, J.-F. (2013).
\newblock Some limit theorems for hawkes processes and application to financial
  statistics.
\newblock \emph{Stochastic Processes and their Applications} \textbf{123},
  2475--2499.

\bibitem[{Bacry \emph{et~al.}(2015)Bacry, Mastromatteo \& Muzy}]{BACRYFINANCE}
Bacry, E., Mastromatteo, I. \& Muzy, J.-F. (2015).
\newblock Hawkes processes in finance.
\newblock \emph{Market Microstructure and Liquidity} \textbf{1}, 1550005.

\bibitem[{Bacry \& Muzy(2014)}]{BACRY}
Bacry, E. \& Muzy, J.-F. (2014).
\newblock Second order statistics characterization of hawkes processes and
  non-parametric estimation.
\newblock \emph{arXiv preprint arXiv:1401.0903} .

\bibitem[{Baraud(2002)}]{baraud}
Baraud, Y. (2002).
\newblock Model selection for regression on a random design.
\newblock \emph{ESAIM: Probability and Statistics} \textbf{6}, 127--146.

\bibitem[{Baraud \emph{et~al.}(2001{\natexlab{a}})Baraud, Comte \&
  Viennet}]{BCV}
Baraud, Y., Comte, F. \& Viennet, G. (2001{\natexlab{a}}).
\newblock Adaptive estimation in autoregression or -mixing regression via model
  selection.
\newblock \emph{Ann. Statist.} \textbf{29}, 839--875.
\newblock \doi{10.1214/aos/1009210692}.

\bibitem[{Baraud \emph{et~al.}(2001{\natexlab{b}})Baraud, Comte \&
  Viennet}]{BCV2}
Baraud, Y., Comte, F. \& Viennet, G. (2001{\natexlab{b}}).
\newblock Model selection for (auto-)regression with dependent data.
\newblock \emph{ESAIM: Probability and Statistics} \textbf{5}, 33–49.

\bibitem[{Bibby \& S{\o}rensen(1995)}]{bibby}
Bibby, B.~M. \& S{\o}rensen, M. (1995).
\newblock Martingale estimation functions for discretely observed diffusion
  processes.
\newblock \emph{Bernoulli} pp. 17--39.

\bibitem[{Bonnet \emph{et~al.}(2018)Bonnet, Rivoirard \& Picard}]{bonnet}
Bonnet, A., Rivoirard, V. \& Picard, F. (2018).
\newblock Modeling spatial genomic interactions with the hawkes model.
\newblock \emph{bioRxiv} p. 214874.

\bibitem[{Br{\'e}maud \& Massouli{\'e}(1996)}]{BM1996}
Br{\'e}maud, P. \& Massouli{\'e}, L. (1996).
\newblock Stability of nonlinear hawkes processes.
\newblock \emph{The Annals of Probability} pp. 1563--1588.

\bibitem[{Carmona \emph{et~al.}(2013)Carmona, Fouque \& Sun}]{CF13}
Carmona, R., Fouque, J.-P. \& Sun, L.-H. (2013).
\newblock Mean field games and systemic risk.
\newblock \emph{Available SSRN} .

\bibitem[{Comte \& Genon-Catalot(2019)}]{CGC2019}
Comte, F. \& Genon-Catalot, V. (2019).
\newblock Nonparametric drift estimation for iid paths of stochastic
  differential equations.
\newblock \emph{hal-02083474} .

\bibitem[{Comte \emph{et~al.}(2007)Comte, Genon-Catalot \& Rozenholc}]{CGCR}
Comte, F., Genon-Catalot, V. \& Rozenholc, Y. (2007).
\newblock Penalized nonparametric mean square estimation of the coefficients of
  diffusion processes.
\newblock \emph{Bernoulli} \textbf{13}, 514--543.

\bibitem[{Daley \& Vere-Jones(2007)}]{DVJ}
Daley, D. \& Vere-Jones, D. (2007).
\newblock \emph{An introduction to the theory of point processes: volume II:
  general theory and structure}.
\newblock Springer Science \& Business Media.

\bibitem[{Dassios \& Zhao(2013)}]{DASSIOS}
Dassios, A. \& Zhao, H. (2013).
\newblock Exact simulation of hawkes process with exponentially decaying
  intensity.
\newblock \emph{Electronic Communications in Probability} \textbf{18}, 1--13.

\bibitem[{Delattre \emph{et~al.}(2016)Delattre, Fournier \& Hoffmann}]{DFH}
Delattre, S., Fournier, N. \& Hoffmann, M. (2016).
\newblock Hawkes processes on large networks.
\newblock \emph{The Annals of Applied Probability} \textbf{26}, 216--261.

\bibitem[{DeVore \& Lorentz(1993)}]{DL93}
DeVore, R. \& Lorentz, G. (1993).
\newblock \emph{Constructive approximation}, vol. 303.
\newblock Springer Science \& Business Media.

\bibitem[{Dion(2014)}]{DION2014}
Dion, C. (2014).
\newblock New adaptive strategies for nonparametric estimation in linear mixed
  models.
\newblock \emph{Journal of Statistical Planning and Inference} \textbf{150}, 30
  -- 48.

\bibitem[{Dion \emph{et~al.}(2019)Dion, Lemler \& L{\"o}cherbach}]{DLL}
Dion, C., Lemler, S. \& L{\"o}cherbach, E. (2019).
\newblock Exponential ergodicity for diffusions with jumps driven by a hawkes
  process.
\newblock \emph{arXiv preprint arXiv:1904.06051} .

\bibitem[{Ditlevsen \& L{\"o}cherbach(2017)}]{DL2016}
Ditlevsen, S. \& L{\"o}cherbach, E. (2017).
\newblock Multi-class oscillating systems of interacting neurons.
\newblock \emph{Stochastic Processes and their Applications} \textbf{127},
  1840--1869.

\bibitem[{Donnet \& Samson(2013)}]{DS2013}
Donnet, S. \& Samson, A. (2013).
\newblock A review on estimation of stochastic differential equations for
  pharmacokinetic/pharmacodynamic models.
\newblock \emph{Advanced Drug Delivery Reviews} \textbf{65}, 929--939.

\bibitem[{Duarte \emph{et~al.}(2016)Duarte, Galves, L{\"o}cherbach \&
  Ost}]{DGL2016}
Duarte, A., Galves, A., L{\"o}cherbach, E. \& Ost, G. (2016).
\newblock Estimating the interaction graph of stochastic neural dynamics.
\newblock \emph{arXiv preprint arXiv:1604.00419} .

\bibitem[{El~Karoui \emph{et~al.}(1997)El~Karoui, Peng \& Quenez}]{ELK}
El~Karoui, N., Peng, S. \& Quenez, M.-C. (1997).
\newblock Backward stochastic differential equations in finance.
\newblock \emph{Mathematical finance} \textbf{7}, 1--71.

\bibitem[{Gloter(2000)}]{GLOTER}
Gloter, A. (2000).
\newblock Discrete sampling of an integrated diffusion process and parameter
  estimation of the diffusion coefficient.
\newblock \emph{ESAIM: Probability and Statistics} \textbf{4}, 205--227.

\bibitem[{Gobet(2002)}]{GOBETIHP}
Gobet, E. (2002).
\newblock Lan property for ergodic diffusions with discrete observations.
\newblock \emph{Ann. Inst. H. Poincar{\'e} Probab. Statist.} \textbf{38},
  711--737.

\bibitem[{Gobet \emph{et~al.}(2004)Gobet, Hoffmann \& Rei{\ss}}]{gobetetal}
Gobet, E., Hoffmann, M. \& Rei{\ss}, M. (2004).
\newblock Nonparametric estimation of scalar diffusions based on low frequency
  data.
\newblock \emph{The Annals of Statistics} \textbf{32}, 2223--2253.

\bibitem[{Gobet \& Matulewicz(2016)}]{GOBET}
Gobet, E. \& Matulewicz, G. (2016).
\newblock Parameter estimation of ornstein--uhlenbeck process generating a
  stochastic graph.
\newblock \emph{Statistical Inference for Stochastic Processes} pp. 1--25.

\bibitem[{Gy{\"o}rfi \emph{et~al.}(2006)Gy{\"o}rfi, Kohler, Krzyzak \&
  Walk}]{GYORFI}
Gy{\"o}rfi, L., Kohler, M., Krzyzak, A. \& Walk, H. (2006).
\newblock \emph{A distribution-free theory of nonparametric regression}.
\newblock Springer Science \& Business Media.

\bibitem[{Hansen \emph{et~al.}(2015)Hansen, Reynaud-Bouret \&
  Rivoirard}]{Hansen15}
Hansen, N., Reynaud-Bouret, P. \& Rivoirard, V. (2015).
\newblock Lasso and probabilistic inequalities for multivariate point
  processes.
\newblock \emph{Bernoulli} \textbf{21}, 83--143.
\newblock \doi{10.3150/13-BEJ562}.

\bibitem[{Has'minskii(1980)}]{Has}
Has'minskii, R. (1980).
\newblock \emph{Stochastic stability of differential equations}.
\newblock Sijthoff \& Noordhoff.

\bibitem[{Hawkes(1971)}]{HAWKES71}
Hawkes, A. (1971).
\newblock Spectra of some self-exciting and mutually exciting point processes.
\newblock \emph{Biometrika} \textbf{58}, 83--90.
\newblock ISSN 00063444.
\newblock \urlprefix\url{http://www.jstor.org/stable/2334319}.

\bibitem[{Hawkes \& Oakes(1974)}]{HO1974}
Hawkes, A. \& Oakes, D. (1974).
\newblock A cluster process representation of a self-exciting process.
\newblock \emph{Journal of Applied Probability} \textbf{11}, 493--503.

\bibitem[{Hoffmann(1999)}]{H99}
Hoffmann, M. (1999).
\newblock Adaptive estimation in diffusion processes.
\newblock \emph{Stochastic processes and their Applications} \textbf{79},
  135--163.

\bibitem[{H\"{o}pfner(2007)}]{HOP}
H\"{o}pfner, R. (2007).
\newblock On a set of data for the membrane potential in a neuron.
\newblock \emph{Mathematical Biosciences} \textbf{207}, 275 -- 301.

\bibitem[{H{\"o}pfner \emph{et~al.}(2016)H{\"o}pfner, L{\"o}cherbach \&
  Thieullen}]{HLT}
H{\"o}pfner, R., L{\"o}cherbach, E. \& Thieullen, M. (2016).
\newblock Ergodicity and limit theorems for degenerate diffusions with time
  periodic drift. application to a stochastic hodgkin- huxley model.
\newblock \emph{ESAIM: Probability and Statistics} \textbf{20}, 527--554.

\bibitem[{Jahn \emph{et~al.}(2011)Jahn, Berg, Hounsgaard \& Ditlevsen}]{JBHD}
Jahn, P., Berg, R., Hounsgaard, J. \& Ditlevsen, S. (2011).
\newblock Motoneuron membrane potentials follow a time inhomogeneous jump
  diffusion process.
\newblock \emph{J. Comput. Neurosci.} \textbf{31}, 563--579.

\bibitem[{Jaisson \& Rosenbaum(2015)}]{RJ2015}
Jaisson, T. \& Rosenbaum, M. (2015).
\newblock Limit theorems for nearly unstable hawkes processes.
\newblock \emph{The Annals of Applied Probability} \textbf{25}, 600--631.

\bibitem[{Kessler \& S{\o}rensen(1999)}]{kesslersorensen}
Kessler, M. \& S{\o}rensen, M. (1999).
\newblock Estimating equations based on eigenfunctions for a discretely
  observed diffusion process.
\newblock \emph{Bernoulli} \textbf{5}, 299--314.

\bibitem[{Kirchner(2017)}]{kirchner2017estimation}
Kirchner, M. (2017).
\newblock An estimation procedure for the hawkes process.
\newblock \emph{Quantitative Finance} \textbf{17}, 571--595.

\bibitem[{Klein \& Rio(2005)}]{klein}
Klein, T. \& Rio, E. (2005).
\newblock Concentration around the mean for maxima of empirical processes.
\newblock \emph{The Annals of Probability} \textbf{33}, 1060--1077.

\bibitem[{Le~Gall(2010)}]{LEGALL}
Le~Gall, J. (2010).
\newblock Calcul stochastique et processus de markov.
\newblock \emph{Notes de cours} .

\bibitem[{Lemonnier \& Vayatis(2014)}]{LV}
Lemonnier, R. \& Vayatis, N. (2014).
\newblock Nonparametric markovian learning of triggering kernels for mutually
  exciting and mutually inhibiting multivariate hawkes processes.
\newblock In \emph{Joint European Conference on Machine Learning and Knowledge
  Discovery in Databases}, pp. 161--176. Springer.

\bibitem[{Lukasik \emph{et~al.}(2016)Lukasik, Srijith, Vu, Bontcheva, Zubiaga
  \& Cohn}]{twit}
Lukasik, M., Srijith, P., Vu, D., Bontcheva, K., Zubiaga, A. \& Cohn, T.
  (2016).
\newblock Hawkes processes for continuous time sequence classification: an
  application to rumour stance classification in twitter.
\newblock \emph{Proceedings of the 54th Annual Meeting of the Association for
  Computational Linguistics (Volume 2: Short Papers)} \textbf{2}, 393--398.

\bibitem[{Mancini(2009)}]{mancini}
Mancini, C. (2009).
\newblock Non-parametric threshold estimation for models with stochastic
  diffusion coefficient and jumps.
\newblock \emph{Scandinavian Journal of Statistics} \textbf{36}, 270--296.

\bibitem[{Masuda(2007)}]{MASUDA}
Masuda, H. (2007).
\newblock Ergodicity and exponential $\beta$-mixing bounds for multidimensional
  diffusions with jumps.
\newblock \emph{Stochastic processes and their applications} \textbf{117},
  35--56.

\bibitem[{Parisi \& Sourlas(1992)}]{PHY}
Parisi, G. \& Sourlas, N. (1992).
\newblock \emph{Supersymmetric field theories and stochastic differential
  equations}.
\newblock World Scientific.

\bibitem[{Rambaldi \emph{et~al.}(2015)Rambaldi, Pennesi \& Lillo}]{RAMBALDI}
Rambaldi, M., Pennesi, P. \& Lillo, F. (2015).
\newblock Modeling foreign exchange market activity around macroeconomic news:
  Hawkes-process approach.
\newblock \emph{Phys. Rev. E} \textbf{91}, 012819.

\bibitem[{Renault \& Touzi(1996)}]{TOUZI}
Renault, E. \& Touzi, N. (1996).
\newblock Option hedging and implied volatilities in a stochastic volatility
  model 1.
\newblock \emph{Mathematical Finance} \textbf{6}, 279--302.

\bibitem[{Reynaud-Bouret \emph{et~al.}(2013)Reynaud-Bouret, Rivoirard \&
  Tuleau-Malot}]{RBRTM}
Reynaud-Bouret, P., Rivoirard, V. \& Tuleau-Malot, C. (2013).
\newblock Inference of functional connectivity in neurosciences via hawkes
  processes.
\newblock In \emph{Global Conference on Signal and Information Processing
  (GlobalSIP), 2013 IEEE}, pp. 317--320. IEEE.

\bibitem[{Reynaud-Bouret \& Roy(2007)}]{RRB}
Reynaud-Bouret, P. \& Roy, E. (2007).
\newblock Some non asymptotic tail estimates for hawkes processes.
\newblock \emph{Bulletin of the Belgian Mathematical Society-Simon Stevin}
  \textbf{13}, 883--896.

\bibitem[{Schmisser(2014{\natexlab{a}})}]{SCHMISSER3}
Schmisser, E. (2014{\natexlab{a}}).
\newblock Nonparametric adaptive estimation of the drift for a jump diffusion
  process.
\newblock \emph{Stochastic Processes and ther Applications} pp. 883--914.

\bibitem[{Schmisser(2014{\natexlab{b}})}]{SCHMISSER2}
Schmisser, E. (2014{\natexlab{b}}).
\newblock Nonparametric estimation of coefficients of a diffusion with jumps.
\newblock \emph{Hal} .

\bibitem[{Shimizu \& Yoshida(2006)}]{SY}
Shimizu, Y. \& Yoshida, N. (2006).
\newblock Estimation of parameters for diffusion processes with jumps from
  discrete observations.
\newblock \emph{Statistical Inference for Stochastic Processes} \textbf{9},
  227--277.

\bibitem[{Tankov \& Voltchkova(2009)}]{TANKOV}
Tankov, P. \& Voltchkova, E. (2009).
\newblock Jump-diffusion models: a practitioner’s guide.
\newblock \emph{Banque et March{\'e}s} \textbf{99}, 24.

\bibitem[{Vasicek(1977)}]{vasicek}
Vasicek, O. (1977).
\newblock An equilibrium characterization of the term structure.
\newblock \emph{Journal of financial economics} \textbf{5}, 177--188.

\bibitem[{Vere-Jones \& Ozaki(1982)}]{Vere82}
Vere-Jones, D. \& Ozaki, T. (1982).
\newblock Some examples of statistical estimation applied to earthquake data.
\newblock \emph{Annals of the Institute of Statistical Mathematics}
  \textbf{34}, 189--207.

\bibitem[{Veretennikov(1997)}]{veret}
Veretennikov, A. (1997).
\newblock On polynomial mixing bounds for stochastic differential equations.
\newblock \emph{Stochastic Processes and their Applications} \textbf{70}, 115
  -- 127.
\newblock ISSN 0304-4149.

\end{thebibliography}

\end{document}